\input ./style/arxiv-general.cfg
\documentclass[aop,MSNbibl,dvips]{arximspdf}
\makeatletter
   \@ifpackageloaded{graphicx}{}{\usepackage{graphicx}}
\makeatother

\doi{10.1214/14-AOP994}
\volume{44}
\issue{2}
\pubyear{2016}
\firstpage{984}
\lastpage{1012}
\docsubty{FLA}

\makeatletter
\def\ssum{\mathop{\sum\sum}}
\newproclaim{ass}{Assumption}

\newcommand{\rrvert}{\vert}

\newcommand{\llvert}{\vert}
\newcommand{\eqref}[1]{(\ref{#1})}

\newtheorem{theo}{Theorem}[section]

\newproclaim{de}[theo]{Definition}
\newtheorem{lem}[theo]{Lemma}
\newproclaim{exxa}[theo]{Example}
\newproclaim{example}[theo]{Example}

\newproclaim{remark}[theo]{Remark}
\newproclaim{remarks}[theo]{Remarks}
\newtheorem{prop}[theo]{Proposition}

\makeatother

\begin{document}
\begin{frontmatter}

\title{Extremes of a class of nonhomogeneous Gaussian~random~fields}
\runtitle{Extremes of Gaussian fields}

\begin{aug}
\author[A]{\fnms{Krzysztof}~\snm{D\c{e}bicki}\thanksref{T1,T2,m1}\ead[label=e1]{Krzysztof.Debicki@math.uni.wroc.pl}},
\author[B]{\fnms{Enkelejd}~\snm{Hashorva}\thanksref{T1,m2}\ead[label=e2]{Enkelejd.Hashorva@unil.ch}}
\and
\author[B]{\fnms{Lanpeng}~\snm{Ji}\corref{}\thanksref{T1,m2}\ead[label=e3]{Lanpeng.Ji@unil.ch}}
\thankstext{T1}{Supported by SNSF Grant 200021-140633/1.}
\thankstext{T2}{Supported by NCN Grant No. 2013/09/B/ST1/01778 (2014-2016).}
\runauthor{K. D{\normalfont{\c{E}}}bicki, E. Hashorva and L. Ji}
\affiliation{University of Wroc\l aw\thanksmark{m1} and University of
Lausanne\thanksmark{m2}}
\address[A]{K. D\c{e}bicki\\
Mathematical Institute\\
University of Wroc\l aw \\
pl. Grunwaldzki 2/4, 50-384 Wroc\l aw\\
Poland\\
\printead{e1}}

\address[B]{E. Hashorva\\
L. Ji\\
Department of Actuarial Science\\
University of Lausanne\\
UNIL-Dorigny, 1015 Lausanne\\
Switzerland\\
\printead{e2}\\
\phantom{E-mail:\ }\printead*{e3}}
\end{aug}

%
\received{\smonth{8} \syear{2013}}
%
\revised{\smonth{10} \syear{2014}}

%
\begin{abstract}
This contribution establishes exact tail
asymptotics of $\sup_{(s,t) \in\mathbf{E}}$ $X(s,t)$ for a large class
of nonhomogeneous Gaussian random fields $X$ on a bounded convex
set $\mathbf{E}\subset\mathbb{R}^2$, with variance function that
attains its
maximum on a segment on $\mathbf{E}$. These findings extend the
classical results for homogeneous Gaussian random fields and
Gaussian random fields with unique maximum point of the variance.
Applications of our result include the derivation of the {exact}
tail asymptotics of the Shepp statistics for stationary Gaussian
processes, Brownian bridge and fractional Brownian motion as well
as the {exact} tail asymptotic expansion
for the maximum loss and span of stationary Gaussian processes.
\end{abstract}

%

%
\end{frontmatter}

\section{Introduction}\label{sec1}
Consider the fractional Brownian motion (fBm) incremental random field
\[
X_\alpha(s,t)=B_{\alpha}(s+t)-B_{\alpha}(s)
,\qquad
(s,t)\in[0,\infty)^2,
\]
where $\{B_\alpha(t),t\in\mathbb{R}
\}$ is a standard fBm with Hurst index $\alpha/2 \in(0,1]$
which is a centered self-similar Gaussian process with stationary
increments and covariance function
\[
\operatorname{Cov}\bigl(B_\alpha(t),B_\alpha(s)\bigr)=
\tfrac{1}{2}\bigl(\llvert t \rrvert ^{\alpha
}+\llvert s \rrvert
^{\alpha}-\mid t-s\mid^{\alpha}\bigr),\qquad s,t\in\mathbb{R}.
\]
For the case $\alpha=1$,
both $X_\alpha(s,t)$ and its standardized version $X_\alpha
^*(s,t)=X_\alpha(s,t)
/t^{\alpha/2}$ appear naturally as limit models; see, for example,
\cite{Cressie}. In the literature,
\[
Y_\alpha(t)=\sup_{s\in[0,S]} X_\alpha(s,t)
\]
is referred to as the Shepp statistics of fBm, whereas $Y_\alpha
^*(t)=\sup_{s\in[0,S]} X_\alpha^*(s,\break t)$ as the standardized
Shepp statistics. Distributional results for $Y^*_1$ are derived in~\cite{Shepp71};
see also \cite{Shepp66} and Theorem~3.2 in \cite{Cressie}. Other
important results for the Shepp statistics
of Brownian motion and related quantities are presented in \cite
{MR905337,MR1833961,MR1331667}.
The first known result for the extremes of the Shepp statistics of
Brownian motion
goes back to \cite{Zholud}, which is {complemented} in
\cite
{Hash13} for the case of fBm
with $\alpha\in(0,1)$. In view of the aforementioned papers for any
$\alpha\in(0,1]$,
%
\begin{equation}
\label{ZH} \mathbb{P} \Bigl( \sup_{(s,t)\in[0,1]^2}
X_\alpha(s,t)> u \Bigr) = C_\alpha u^{4/\alpha-2} \Psi(u)
\bigl(1+o(1)\bigr)
\end{equation}
holds as $u\to\infty$ with
$C_\alpha$
a positive constant
and $\Psi(\cdot)$ the survival function of an $N(0,1)$ random variable.
There is no result for the case $\alpha\in(1,2)$ in the
literature; we shall cover this gap in Proposition~\ref{propfBm}.

Results for the tail asymptotics of supremum of the standardized Shepp
statistics can be derived using the findings of
\cite{ChanLai} and \cite{MikhalevaPit}; see also \cite{KabB,KabA}.
However, 
this is not the case for the tail asymptotics of the supremum of the
Shepp statistics $Y_\alpha$;
no theoretical results in the literature can be applied for this case.
This is due to the fact that on $[0,1]^2$ the
variance of $X_\alpha$ attains its maximum at an infinite number
of points, that is,
its maximal value is attained for any $s\in[0,1]$ and $t=1$.

In the asymptotic theory of Gaussian random fields,
{if} the random field has a nonconstant variance function, which
attains its maximum at a unique (or finite) number of points, then
under the so-called Piterbarg conditions, the exact tail asymptotics of
supremum of Gaussian random fields with certain $(E,\alpha)$ structures
for the variance and the correlation functions are derived by
relying on the Double--Sum method; see, for example, the standard
monograph \cite{Pit96}.


The principle aim of this contribution is to extend Piterbarg's
asymptotic theory for Gaussian random fields
to the case where the maximum of the variance function on a bounded
convex set $\mathbf{E}$ is attained on finite number of disjoint
segments on $\mathbf{E}$.
In particular, we assume that $\{X(s,t), (s,t) \in\mathbf{E}\}$,
$\mathbf{E}=[0,S] \times[0,T], S,T>0$, is a centered Gaussian random
field with
variance function $\sigma^2(s,t)=\operatorname{Var}(X(s,t))$ that
satisfies the
following assumption.

\renewcommand{\theass}{A\arabic{ass}}
\begin{ass}\label{assA1}
There exists some positive function $\sigma(t)$
which attains its unique maximum on $[0,T]$ at $T$, and further
%
\begin{eqnarray}
\label{var} \sigma(s,t)&=&\sigma(t)\qquad \forall(s,t)\in\mathbf{E},
\nonumber
\\[-8pt]
\\[-8pt]
\nonumber
 \sigma(t)&=& 1- b(T-t) ^\beta\bigl(1+o(1)\bigr),\qquad t\uparrow T
\end{eqnarray}
hold for some $\beta,b >0$.
\end{ass}

We shall impose the following assumption on the correlation
function $r(s,t, s',t')=\mathbb{E} (\overline{X}(s,t)\overline
{X}(s',t') )$
where $\overline{X}(s,t)={X(s,t)}/{\sigma(s,t)}$:

\begin{ass}\label{assA2}
There exist constants $a_1>0,\break a_2>0,a_3\neq0$ and
$\alpha_1, \alpha_2 \in(0,2]$ such that
%
\begin{eqnarray}
\label{eq:r}&& r\bigl(s,t,s',t'\bigr)
\nonumber
\\[-8pt]
\\[-8pt]
\nonumber
&&\qquad= 1- \bigl(
\bigl\llvert a_1\bigl(s-s'\bigr) \bigr\rrvert
^{\alpha
_1}
+\bigl\llvert a_2\bigl(t-t'\bigr)+a_3
\bigl(s-s'\bigr) \bigr\rrvert ^{\alpha_2} \bigr) \bigl(1+ o(1)
\bigr)
\end{eqnarray}
holds uniformly with respect to $s,s'\in[0,S]$, as $\llvert  s-s'
\rrvert \rightarrow0,
t,t' \uparrow T$, and further, there exists some constant $\delta_0\in
(0,T)$ such that
%
\begin{equation}
\label{eq:r1} r\bigl(s,t,s',t'\bigr)< 1
\end{equation}
holds for any $s,s'\in[0,S]$ satisfying $s\neq s'$, and $
t,t'\in[\delta_0,T]$. 
\end{ass}

Note that in \ref{assA2} we assume that $a_3\neq0$, which includes a
large class of correlation functions with $(E, \alpha)$ structure
dealt with in \cite{Pit96}; the classical case $a_3=0$ is discussed in
Remark~\ref{remarka3}.

Our main result, presented in {Theorem~\ref{thM}} (and stated in higher
generality in Remarks \ref{remarkmain}),
derives the exact tail asymptotic behavior of supremum of
nonhomogeneous Gaussian random
fields $X$ satisfying \ref{assA1} and \ref{assA2} and a H\"{o}lder condition
formulated below in Assumption \ref{assA3}.
As an illustration to the derived theory, we analyze exact asymptotics
of the tail distribution of extremes of Shepp statistics, the maximum
loss and the span for a large class of Gaussian processes.

Organization of the paper: Our principal
findings are presented in Section~\ref{sec2} followed
by two sections dedicated to applications and examples. All the proofs
are relegated to Section~\ref{sec5} and the \hyperref[app]{Appendix}.

\section{Main results}\label{sec2}

In this section, we are concerned about the asymptotics of
\[
\mathbb{P} \Bigl( \sup_{(s, t)\in\mathbf{E}}X(s,t)>u \Bigr) ,\qquad u\rightarrow
\infty
\]
discussing first the case that $\mathbf{E}=[0,S]\times[0,T]$.

The Pickands and Piterbarg lemmas (cf. \cite{Pit96}) are
fundamental in the analysis of the tail asymptotic
behavior of
supremum of nonsmooth centered Gaussian processes and Gaussian random
fields. Restricting ourselves to the case that $\{X(t),t\ge0\}
$ is a centered stationary Gaussian process
with a.s. continuous sample paths and correlation function $r(t)$,
such that $r(t) = 1 - t^{\alpha}(1+o(1))$ as $t \rightarrow0$, with
$\alpha\in(0, 2]$, and $r(t)<1$ for all $t>0$,
in\break view of the seminal papers by J. Pickands III (see \cite{PicandsB,PicandsA}), for any $T\in(0,\infty)$
%
\begin{equation}
\label{pic} \mathbb{P} \Bigl( \sup_{t\in[0,T]}X(t)>u \Bigr) =
\mathcal {H}_{\alpha}T u^{{2}/{\alpha}}\Psi(u) \bigl(1+o(1)\bigr) ,\qquad u
\rightarrow\infty.
\end{equation}
Here, $ \mathcal{H}_{\alpha}$ is the \textit{Pickands constant} defined by
\[
\label{pick} \mathcal{H}_{\alpha}=\lim_{T\rightarrow\infty}
\frac{1}{T} \mathcal{H}_{\alpha}[0, T] \in(0,\infty)
\]
with
\[
\mathcal{H}_{\alpha}[0, T]=\mathbb{E} \Bigl(\exp \Bigl(\sup
_{t\in[0,T]} \bigl(\sqrt{2}B_\alpha(t)-t^{\alpha}
\bigr) \Bigr) \Bigr).
\]
The derivation of \eqref{pic} is based on \textit{Pickand's lemma}
which states that
%
\begin{equation}
\label{eq:pick1} \mathbb{P} \Bigl( \sup_{t\in[0,u^{-{2}/{\alpha}}T]}X(t)>u \Bigr) =
\mathcal{H}_{\alpha}[0, T] \Psi(u) 
\bigl(1+o(1)\bigr),\qquad  u
\rightarrow\infty.
\end{equation}
In \cite{Pit72}, Piterbarg rigorously proved Pickand's theorem and
further derived
a crucial extension of \eqref{eq:pick1} which
we shall refer to as the \textit{Piterbarg lemma}; it states that
%
\begin{equation}
\label{eq:piter1} \mathbb{P} \biggl( \sup_{t\in[0,u^{-{2}/{\alpha}}T]}\frac
{X(t)}{1+bt^\alpha}>u
\biggr) =\mathcal{P}_{\alpha}^b[0,T]\Psi(u)%
\bigl(1+o(1)\bigr),\qquad u \rightarrow\infty
\end{equation}
holds for any $b>0$ with
\[
\label{pick} \mathcal{P}_{\alpha}^b[0, T]=\mathbb{E} \Bigl(
\exp \Bigl(\sup_{t\in[0,T]} \bigl(\sqrt{2}B_\alpha(t)-(1+b)t^{\alpha}
\bigr) \Bigr) \Bigr) \in(0,\infty).
\]
The positive constant (referred to as the \textit{Piterbarg constant}) given by
\[
\mathcal{P}_\alpha^b=\lim_{T\rightarrow\infty}\mathcal
{P}_{\alpha}^b[0, T] \in (0,\infty)
\]
appears naturally when dealing with the extremes of nonstationary
Gaussian processes or Gaussian random fields; see, for example,
\cite{Pit96} and our main result below. It is known that $\mathcal{H}_1=1$,
$\mathcal{H}_2={1}/{\sqrt{\pi}}$, and
%
\begin{equation}
\label{eqpp} \mathcal{P}_1^b=1+\frac{1}{b},\qquad
\mathcal{P}_2^b=\frac
{1}{2} \biggl(1+\sqrt{1+
\frac{1}{b}} \biggr),\qquad b>0
\end{equation}
see, for example, \cite{Albin1990,MR1993262,dieker2005extremes,MR3091101,Harper2,DikerY}.

We note in passing that for stationary Gaussian processes \cite
{AlbinC} and \cite{Berman82} presented new elegant proofs of \eqref{pic}
without using the Pickands lemma.
The following extension of the Pickands and Piterbarg lemmas plays an
important role in our analysis.
Hereafter, {we} denote by $\tilde{B}_{ \alpha}$ and $B_{\alpha}$
two independent fBm's defined on $\mathbb{R}$ with Hurst index $\alpha
/2\in(0,1]$.
Recall that $\Psi(\cdot)$ denotes the survival function of an
$N(0,1)$ random variable; we write below $\Gamma(\cdot)$ for the
Euler Gamma function.

\begin{lem}\label{Lem2}
Let $\{\eta(s,t), (s,t)\in[0,\infty)^2\}$ be a centered
homogeneous
Gaussian random field with covariance function
\[
r_\eta(s,t)=\exp \bigl(-\llvert a_1s \rrvert
^{\alpha
_1}-\llvert a_2 t-a_3s \rrvert
^{\alpha
_2} \bigr), \qquad (s,t)\in[0,\infty)^2,
\]
where constants $ \alpha_i \in(0,2], i=1,2$, $a_1>0,a_2>0, a_3\in
\mathbb{R}$.
Let further $b, S,T$ be three positive constants. If $\beta\ge\alpha
_2\ge\alpha_1$, then for any positive measurable function $g(u),
u>0$ satisfying $\lim_{u\rightarrow\infty}g(u)/u=1$
%
\begin{eqnarray}
\label{eq:lem1} &&\mathbb{P} \biggl( \sup_{(s,t)\in[0,Su^{-{2}/{\alpha
_1}}]\times[0,Tu^{-{2}/{\alpha_2}}]}
\frac{\eta
(s,t)}{1+bt^{\beta}}>g(u) \biggr)
\nonumber
\\[-8pt]
\\[-8pt]
\nonumber
&&\qquad =\mathcal{H}_{Y}^b[S,T]\Psi\bigl(g(u)\bigr)
\bigl(1+o(1)\bigr), \qquad u\rightarrow\infty,
\end{eqnarray}
where
%
\begin{eqnarray}
\label{eq:pick} \mathcal{H}_Y^{d}[S,T]&=&\mathbb{E}
\Bigl(\exp \Bigl(\sup_{(s,t)\in
[0,S]\times
[0,T]} \bigl(\sqrt{2}Y(s,t)
 -\sigma_Y^2(s,t)-d(t) \bigr) \Bigr) \Bigr)
 \nonumber
 \\[-8pt]
 \\[-8pt]
 \nonumber
 &\in&(0,
\infty)
\end{eqnarray}
with $\sigma^2_Y(s,t)=\operatorname{Var}(Y(s,t))$ and
%
\begin{eqnarray}
\label{eqY} Y(s,t)&=& \cases{ %
Y_1(s,t):=
\tilde{B}_{\alpha_1}(a_1 s)+ B_{\alpha_2}(a_2
t-a_3 s), & \quad $\alpha_1=\alpha_2,$
\vspace*{2pt}\cr
Y_2(s,t):= \tilde{B}_{\alpha_1}(a_1
s)+B_{\alpha_2}(a_2 t), & \quad $\alpha _1<
\alpha_2,$} 
\nonumber
\\[-8pt]
\\[-8pt]
\nonumber
d(t)&=& \cases{
0, & \quad $\beta>\alpha_2,$
\vspace*{2pt}\cr
bt^\beta, &\quad $\beta=\alpha_2,$}\qquad
 (s,t)
\in[0,\infty)^2.
\end{eqnarray}
\end{lem}

Using the definition of $Y_1$ and $Y_2$ appearing in \eqref{eqY}
we shall determine, for given $a_i$'s, $\alpha_i$'s and $b, \beta$ as above,
the following constants (referred to as \textit{generalized
Pickands--Piterbarg constants}):
\[
\mathcal{M}_{Y,\beta}^{b}=\lim
_{T\rightarrow\infty}\lim_{
S\rightarrow\infty}\frac{1}{S}
\mathcal{H} _{Y}^{b}[S,T]\in(0,\infty)
\]
and
\begin{eqnarray*}
\widetilde{\mathcal{M}}_{Y,\beta}^{b}&=&\lim
_{T\rightarrow\infty
}\lim_{ S\rightarrow\infty
}\frac
{1}{S} \mathbb{E}
\Bigl(\exp \Bigl(\sup_{(s,t)\in[0,S]\times
[-T,T]} \bigl(\sqrt{2}Y(s,t)
 -\sigma_{Y}^2(s,t)-bt^\beta \bigr) \Bigr)
\Bigr) \\
&\in&(0,\infty). 
\end{eqnarray*}
Here, $\mathcal{M}_{Y,\beta}^{b}$ and $\widetilde{\mathcal
{M}}_{Y,\beta
}^{b}$ are
defined only for $\beta=\alpha_2$. 
Note that we suppress $ a_i$'s and $\alpha_i$'s in the definition of
$\mathcal{M}_{Y,\beta}^{b}$ and $\widetilde{\mathcal{M}}_{Y,\beta}^{b}$
since they appear directly in the definition of $Y$.

Additional to \ref{assA1} and \ref{assA2} we shall impose the following H\"
older condition, which in the literature is called \textit{regularity};
see \cite{Pit96}.

\begin{ass}\label{assA3}
There exist positive constants $\rho_1,\rho_2,
\gamma, \mathcal{Q}$ such that
\[
\mathbb{E} \bigl(\bigl(X(s,t)-X\bigl(s',t'
\bigr)\bigr)^{2} \bigr) \leq \mathcal {Q} \bigl(\bigl|t-t'\bigr|^{\gamma}+
\bigl\llvert s-s' \bigr\rrvert ^{\gamma
} \bigr)
\]
holds for all $ t,t'\in[\rho_1, T], s,s' \in[0,S]$ satisfying $\llvert  s-s' \rrvert <\rho_2$.
\end{ass}

We present next our main result.


\begin{theo}\label{thM} Let $ \{X(s,t), (s,t) \in\mathbf{E}\},
\mathbf{E}=[0,S]
\times
[0,T]$ be a centered
Gaussian random field with a.s. continuous sample paths. Suppose that
Assumptions \ref{assA1}--\ref{assA3} are satisfied
with the parameters mentioned therein. Then, as \mbox{$u\rightarrow\infty$},
\begin{longlist}[(iii)]
\item[(i)] if $\beta>\max(\alpha_1,\alpha_2)$
%
\begin{eqnarray}
\label{main1} \mathbb{P} \Bigl( \sup_{(s,t)\in\mathbf{E}} X(s,t)> u \Bigr)
&=& S \Gamma ({1}/{\beta }+1 ) \prod_{k=1}^2(a_k
\mathcal{H}_{\alpha_k}) b^{-{1}/{\beta
}}
\nonumber
\\[-8pt]
\\[-8pt]
\nonumber
&&{} \times u^{{2}/{\alpha_2}+{2}/{\alpha_1}-{2}/{\beta}} \Psi(u) \bigl(1+o(1)\bigr);
\end{eqnarray}
\item[(ii)] if $\beta=\alpha_2=\alpha_1$
%
\begin{equation}
\mathbb{P} \Bigl( \sup_{(s,t)\in\mathbf{E}} X(s,t)> u \Bigr) = S
\mathcal{M}_{Y_1,\alpha_1}^{b} u^{{2}/{\alpha_1}} 
\Psi(u)
\bigl(1+o(1)\bigr);
\end{equation}
\item[(iii)] if $\beta=\alpha_2>\alpha_1$
%
\begin{equation}
\label{main2} \mathbb{P} \Bigl( \sup_{(s,t)\in\mathbf{E}} X(s,t)> u \Bigr) =
S a_1a_2\mathcal{P}_{\alpha
_2}^{ba_2^{-\alpha_2}}
\mathcal{H}_{\alpha_1} u^{{2}/{\alpha
_1}} 
\Psi(u) \bigl(1+o(1)\bigr);
\end{equation}
\item[(iv)] if $\beta<\alpha_2=\alpha_1$
%
\begin{equation}
\mathbb{P} \Bigl( \sup_{(s,t)\in\mathbf{E}} X(s,t)> u \Bigr) = S
\bigl(a_1^{\alpha
_1}+\llvert a_3 \rrvert
^{\alpha_1}\bigr)^{{1}/{\alpha_1}} \mathcal{H}_{\alpha_1} u^{
{2}/{\alpha
_1}}
\Psi(u) \bigl(1+o(1)\bigr);
\end{equation}
\item[(v)] if $\beta< \alpha_2 $ and $\alpha_1<\alpha_2$
%
\begin{equation}
\label{main3} \mathbb{P} \Bigl( \sup_{(s,t)\in\mathbf{E}} X(s,t)> u \Bigr) =
S a_1 \mathcal{H}_{\alpha_1} u^{
{2}/{\alpha
_1}} 
\Psi(u)
\bigl(1+o(1)\bigr);
\end{equation}
\item[(vi)] if $\beta=\alpha_1>\alpha_2$
%
\begin{equation}\quad
\mathbb{P} \Bigl( \sup_{(s,t)\in\mathbf{E}} X(s,t)> u \Bigr) = S
a_1\mathcal{P}_{\alpha
_1}^{b
({\llvert  a_3 \rrvert }/{(a_1a_2)} )^{\alpha_1}} \mathcal{H}_{\alpha_2}
u^{
{2}/{\alpha
_2}} 
\Psi(u) \bigl(1+o(1)\bigr);
\end{equation}
\item[(vii)] if $\beta< \alpha_1 $ and $\alpha_2<\alpha_1$
\[
\mathbb{P} \Bigl( \sup_{(s,t)\in\mathbf{E}} X(s,t)> u \Bigr) = S \llvert
a_3 \rrvert \mathcal{H}_{\alpha_2} u^{
{2}/{\alpha_2}} 
\Psi(u) \bigl(1+o(1)\bigr).
\]
\end{longlist}
\end{theo}

\begin{remark}\label{remarka3}
If $a_3=0$, then there are only three scenarios to be
considered.
In particular
if $\beta>\alpha_2$, then \eqref{main1} holds. If $\beta=\alpha_2$,
then~\eqref{main2} holds, whereas if $\beta<\alpha_2$,
then \eqref{main3} is valid.
\end{remark}

\begin{remark}\label{remarkmain} (a) Let $\mathbf{E}$ be any
bounded convex
subset of $\mathbb{R}^2$. Assume that on $\mathbf{E}$ the maximum
of the standard
deviation $\sigma(s,t)$ is attained only on a segment $\mathbf
{L}$ which is
inside of $\mathbf{E}$, parallel to $s$-axis and of length $\ell$.
Then the claims of {Theorem~\ref{thM}} are still valid, by replacing
$S$ with
$\ell$ in cases (i)--(vii), $\Gamma(\cdot)$ with $2\Gamma(\cdot)$ in
case (i), $\mathcal{M}_{Y_1,\alpha_1}^{b}$ with $\widetilde
{\mathcal{M}
}_{Y_1,\alpha_1}^{b}$ in cases (ii), $\mathcal{P}_{\alpha
_2}^{ba_2^{-\alpha
_2}} $ with $\widetilde{\mathcal{P}}_{\alpha_2}^{ba_2^{-\alpha_2}}
$ in case
(iii), and $\mathcal{P}_{\alpha_1}^{b (\llvert  a_3 \rrvert /(a_1a_2) )^{\alpha_1}}$
with $\widetilde{\mathcal{P}}_{\alpha_1}^{b ( \llvert  a_3
\rrvert /(a_1a_2)
)^{\alpha
_1}}$ in case (vi), respectively. Here, $\widetilde{\mathcal
{P}}_{\alpha
}^{b}$, with $b>0$ and $\alpha\in(0,2]$ is the Piterbarg constant
defined on the real line, that is,
\[
\widetilde{\mathcal{P}}_{\alpha}^b=\lim_{T\rightarrow\infty
}
\mathbb{E} \Bigl(\exp \Bigl(\sup_{t\in[-T,T]} \bigl(\sqrt
{2}B_\alpha(t)-(1+b)t^{\alpha} \bigr) \Bigr) \Bigr) \in (0,
\infty).
\]

(b) Assume that on $\mathbf{E}$ the maximum of the standard deviation
$\sigma
(s,t)$ is attained only on $n$ segments $\{\mathbf{L}_i\}_{i=1}^
n$ which are
inside or on the boundary of $\mathbf{E}$, and parallel to
$s$-axis. By the
convexity of $\mathbf{E}$, we can always find $n$ nonadjacent
convex sets
$\{\mathbf{E} _i\}_{i=1}^ n$ such that $\mathbf{L}_i\subset
\mathbf{E} _i\subset\mathbf{E}$,
$i=1,\ldots,n$. If further
for any $i\neq j$
%
\begin{equation}
\label{eq:r2'} \sup_{(s,t)\in\mathbf{E} _i,(s',t')\in\mathbf{E}
_j}r\bigl(s,t,s',t'
\bigr)< 1
\end{equation}
holds, then
%
\begin{equation}
\label{eq:EE} \mathbb{P} \Bigl( \sup_{(s,t)\in\mathbf{E} } X(s,t)> u \Bigr) =
\sum_{i=1}^n\mathbb{P} \Bigl( \sup
_{(s,t)\in
\mathbf{E} _i} X(s,t)> u \Bigr) \bigl(1+o(1)\bigr)
\end{equation}
as $u\to\infty$. Additionally, suppose that on each $\{\mathbf
{E} _i\}_{i=1}^ n$
the Assumptions \ref{assA1}--\ref{assA3} are satisfied. Then an explicit
expression for \eqref{eq:EE} can be established by applying the results
in {Theorem~\ref{thM}} and Remark \ref{remarkmain}(a) above.

(c) {Similar} results can also be obtained when the segments $\{
\mathbf{L}
_i\}_{i=1}^ n$, where the maximum of $\sigma(s,t)$ is attained,
are nonparallel and disjoint. Specifically, 
we see from Remark \ref{remarkmain}(b) that it is sufficient to consider the asymptotics of
\[
\mathbb{P} \Bigl( \sup_{(s,t)\in\mathbf{E} _i} X(s,t)> u \Bigr),\qquad u\rightarrow
\infty, i=1,\ldots,n,
\]
respectively.
Let $(s,t)^\top$ be the transpose of $(s,t)$. Then, for any
$i=1,\ldots,n$, there is a nondegenerate lower triangular (rotation)
matrix $A_i\in\mathbb{R}^{2\times2}$ such that the maximum of the
variance of $X((A_i(s,t)^\top)^{\top})$ on $A_i^{-1}\mathbf{E}
_i=\{(\tilde s,\tilde t)\dvtx\break  (\tilde s,\tilde t)^\top=A_i^{-1}(s,t)^\top
, (s,t)\in\mathbf{E} _i\}$ is attained on a line parallel to
$s$-axis or $t$-axis.
Consequently, similar results as in {Theorem~\ref{thM}} can be
obtained if
certain Assumptions as \ref{assA1}--\ref{assA3} are satisfied by each $\{
X((A_i(s,t)^\top)^{\top}), (s,t)\in A_i^{-1}\mathbf{E} _i\}$.
\end{remark}

We conclude this section with an example, which illustrates the
existence of all the cases discussed in {Theorem~\ref{thM}}.

\begin{example} Consider a Gaussian random field defined as
\[
Z(s,t)=\tfrac{1}{\sqrt{2}}\bigl(Y(s+t)-X(s)\bigr) \bigl(1-b(T-t)^\beta
\bigr),\qquad (s,t)\in[0,S]\times[0,T],
\]
where $b,\beta$ are two positive constants, and $X,Y$ are two
independent centered stationary Gaussian processes with covariance
functions $r_X,r_Y$ satisfying as $t\to0$
\[
r_X(t)=1-a_1t^{\alpha_1}\bigl(1+o(1)\bigr),\qquad
r_Y(t)=1-a_2t^{\alpha_2}\bigl(1+o(1)\bigr)
\]
for some constants $a_i>0,\alpha_i\in(0,2], i=1,2$. Further,
assume that
\[
r_X(s)<1, \forall s\in(0,S]\qquad r_Y(t)<1\qquad \forall t
\in(0,S+T].
\]
It follows that the assumptions of {Theorem~\ref{thM}} are satisfied by
$\{Z(s,t), (s,t)\in[0,S]\times[0,T]\}$.
\end{example}

\section{Extremes of Shepp statistics}\label{sec3}
For a given centered Gaussian process $\{X(t),t\ge0\}$, we shall
define the incremental random field $Z$ by
%
\begin{equation}
\label{eq:Z} Z(s,t)=X(s+t)-X(s),\qquad (s,t)\in[0,S]\times[0,T].
\end{equation}
The asymptotic analysis of the supremum of the Shepp statistics
\[
Y(t)= \sup_{s\in[0,S]} Z(s,t),\qquad t\in[0,T]
\]
boils down to the study of the tail asymptotics of the double-supremum
$\sup_{(s,t)\in[0,S]\times[0,T]} Z(s,t)$. In this section, we shall
consider several important examples
which can be analysed utilising the theory developed in Section~\ref{sec2}.

\subsection{Stationary Gaussian processes}\label{sec3.1}
Consider the Gaussian random field $Z$ as in \eqref{eq:Z} where $X$ is
a centered stationary Gaussian process with covariance function
$r_X$ satisfying the following conditions:
\begin{longlist}[S1:]
\item[S1:] $r_X(t)$ attains its minimum on $[0,T]$ at the
unique point $t=T$;

\item[S2:] there exist positive constants $\alpha_1,
a_1,a_2$ and
$\alpha
_2\in(0,2)$ such that 
\begin{eqnarray*}
r_X(t)&=&r_X(T)+a_1(T-t)^{\alpha_1}
\bigl(1+o(1)\bigr),\qquad t\rightarrow T,
\\
 r_X(t)&=&1-a_2t^{\alpha_2}\bigl(1+o(1)\bigr),\qquad t
\rightarrow0;
\end{eqnarray*}

\item[S3:] $r_X(s)<1$ for any $s\in(0,S+T]$. 
\end{longlist}

%
\begin{prop}\label{propSGP} Let $\{Z(s,t), (s,t)\in[0,S]\times[0,T]\}$
be an incremental random field given as in \eqref{eq:Z} with
$r_X$ satisfying \textup{{S1--S3}}.
Suppose that $r_X$ is twice continuously differentiable on $[\mu,T]$
for some $\mu\in(0,T)$,
$\llvert  r^{\prime\prime}_X(T) \rrvert \in(0,\infty)$, and let
$b_i={a_i}/\rho_T^2,i=1,2$ with $\rho_T=\sqrt{2 (1-r_X(T))}$. Then,
as $u\rightarrow\infty$, 
\begin{longlist}
\item[(i)] if $\alpha_1>\alpha_2$
\begin{eqnarray*}
&&\mathbb{P} \Bigl( \sup_{(s,t)\in[0,S] \times[0, T]} Z(s,t)> u \Bigr) \\
&&\qquad= S \Gamma
(1/\alpha _1+1) \mathcal{H}_{\alpha_2}^2
b_2^{{2}/{\alpha_2}} b_1^{-
{1}/{\alpha
_1}} \biggl(\frac{u}{ \rho_T} \biggr)^{{4}/{\alpha
_2}-
{2}/{\alpha_1}} 
\Psi
\biggl(\frac{u}{ \rho_T} \biggr) \bigl(1+o(1)\bigr);
\end{eqnarray*}
\item[(ii)] if $\alpha_1=\alpha_2$
\[
\mathbb{P} \Bigl( \sup_{(s,t)\in[0,S] \times[0, T]} Z(s,t)> u \Bigr) = S \mathcal{M}
_{Y,\alpha
_1}^{b_1} \biggl(\frac{u}{ \rho_T} \biggr)^{{2}/{\alpha_2}}
\Psi \biggl(\frac{u}{ \rho_T} \biggr) \bigl(1+o(1)\bigr),
\]
where
\[
Y(s,t):= \tilde{B}_{\alpha_2} \bigl(b_2^{{1}/{\alpha_2}} s
\bigr)+ B_{\alpha_2} \bigl(b_2^{{1}/{\alpha_2}}
t-b_2^{{1}/{\alpha
_2}} s \bigr), \qquad (s,t)\in[0,\infty)^2;
\]
\item[(iii)] if $\alpha_1<\alpha_2$
\[
\mathbb{P} \Bigl( \sup_{(s,t)\in[0,S] \times[0, T]} Z(s,t)> u \Bigr) = S
(2b_2)^{
{1}/{\alpha_2}} \mathcal{H}_{\alpha_2} \biggl(
\frac{u}{ \rho
_T} \biggr)^{{2}/{\alpha
_2} } 
\Psi \biggl(
\frac{u}{ \rho_T} \biggr) \bigl(1+o(1)\bigr).
\]
\end{longlist}
\end{prop}


We present two important examples that illustrate {Proposition~\ref
{propSGP}}.

\begin{example}[(Slepian process)] Consider $X$ to be the Slepian
process, that is,
\[
X(t)=B_1(t+1)-B_1(t),\qquad t\in[0,\infty),
\]
with $B_1$ the standard Brownian motion. 
It follows that the assumptions of {Proposition~\ref{propSGP}} are satisfied,
hence as $u\to\infty$
\[
\mathbb{P} \Bigl( \sup_{(s,t)\in[0,1] \times[0, {1}/{2}]} Z(s,t)> u \Bigr) = \mathcal{M}
_{Y,1}^1 u^2\Psi(u) \bigl(1+o(1)
\bigr)
\]
holds with $Y(s,t):= \tilde{B}_1 ( s )+ B_1
(t-s ), (s,t)\in
(0,\infty)^2$.
\end{example}

\begin{example}[(Ornstein--Uhlenbeck process)] Consider a centered stationary
Gaussian process $X$ with covariance function
$r(t)=e^{-t},t\ge0$. Then following {Proposition~\ref{propSGP}},
\[
\mathbb{P} \Bigl( \sup_{(s,t)\in[0,1]^2} Z(s,t)> u \Bigr) =
\mathcal{M}_{Y,1}^{b_1}b_1 u^2\Psi (
\sqrt{b_1}u) \bigl(1+o(1)\bigr), \qquad u\rightarrow\infty,
\]
with $b_1=e^{-1}/(2(1-e^{-1}))$, $b_2=
1/(2(1-e^{-1}))$ and
$
Y(s,t):= \tilde{B}_{1} (b_2 s )+
B_{1} (b_2 t-b_2
s
),  (s,t)\in(0,\infty)^2$.
\end{example}

\subsection{Brownian bridge}\label{sec3.2}
In this section, we analyze
%
\begin{equation}
\label{eq:Z2} Z(s,t)=X(s+t)-X(s),\qquad s,s+t\in[0,1],
\end{equation}
where $X(s):=B_1(s)-sB_1(1), s\in[0,1]$ is a Brownian bridge (recall
$B_1$ is a standard Brownian motion). Clearly,
$X$ is nonstationary and, therefore, we cannot apply {Proposition~\ref
{propSGP}}
for this case.

%
\begin{prop}\label{propBB}
If $\{Z(s,t), (s,t)\in[0,1/2]^2\}$ is given by \eqref{eq:Z2}, then
%
\begin{equation}\quad
\mathbb{P} \Bigl( \sup_{(s,t)\in[0,{1}/{2}]^2 } Z(s,t)> u \Bigr) =
2^{
{5}/{2}}\sqrt{\pi} u^3 \Psi(2u) \bigl(1+o(1)\bigr),\qquad u
\rightarrow\infty.
\end{equation}
\end{prop}

\subsection{Fractional Brownian motion}\label{sec3.3}
Consider the fBm incremental random field
%
\begin{equation}
\label{eq:ZH} Z(s,t)=B_\alpha(s+t)-B_\alpha(s),\qquad (s,t)\in[0,S]
\times[0,1],
\end{equation}
where $B_\alpha$ is the fBm with Hurst index $\alpha/2\in(0,1)$.

The following proposition extends the main result of \cite{Hash13} to
the whole range of $\alpha\in(0,2)$.

\begin{prop}\label{propfBm}
Let $\{Z(s,t), (s,t)\in[0,S]\times[0,1]\}$ be given as in
\eqref
{eq:ZH}. We have, as $u\rightarrow\infty$,
\begin{longlist}[(iii)]
\item[(i)] if $\alpha\in(0,1)$
%
\begin{equation}\quad
\mathbb{P} \Bigl( \sup_{(s,t)\in[0,S]\times[0,1]} Z(s,t)> u \Bigr) = S
2^{1-2/\alpha}\alpha^{-1}\mathcal{H}_{\alpha}^2
u^{
{4}/{\alpha}-2 } \Psi(u) \bigl(1+o(1)\bigr);
\end{equation}
\item[(ii)] if $\alpha=1$
%
\begin{equation}
\mathbb{P} \Bigl( \sup_{(s,t)\in[0,S]\times[0,1]} Z(s,t)> u \Bigr) = S \mathcal{M}
_{Y,1}^{
{1}/{2}} u^{ 2 } \Psi(u) \bigl(1+o(1)\bigr),
\end{equation}
with
\[
Y(s,t):= \tilde{B}_{1} \bigl(2^{-1} s \bigr)+
B_{1} \bigl(2^{-1} (t-s) \bigr),\qquad (s,t)\in[0,
\infty)^2;
\]
\item[(iii)] if $\alpha\in(1,2)$
%
\begin{equation}
\mathbb{P} \Bigl( \sup_{(s,t)\in[0,S]\times[0,1]} Z(s,t)> u \Bigr) = S
\mathcal{H}_{\alpha} u^{
{2}/{\alpha}} \Psi(u) \bigl(1+o(1)\bigr).
\end{equation}
\end{longlist}
\end{prop}

\section{Extremes of maximum loss and span of Gaussian processes}\label{sec4}

Let $\{\xi(t),\break t\in[0,1]\}$ be a Gaussian process with a.s.
continuous sample paths.
The maximum loss of the process $\xi$ is given by
\[
\chi_1(\xi)=\max_{0\le s\le t\le1}\bigl( \xi(s)-\xi(t)\bigr),
%
\]
and its span is defined as
\[
\chi_2(\xi)=\max_{t\in[0,1]}\xi(t)-\min
_{t\in[0,1]}\xi (t). 
\]
The notion of the maximum loss of certain Gaussian processes (e.g.,
Brownian motion and fBm, etc.)
plays an important role in finance and insurance modelling; see, for
example, \cite{VardarB}, \cite{VardarA} and references therein.

In this section, as an application of {Theorem~\ref{thM}}
and Remark \ref{remarkmain}, we derive {exact} tail
asymptotics of the maximum loss for both stationary Gaussian process
(in {Proposition~\ref{PP}}) and for
Brownian bridge (in {Proposition~\ref{BB2}}).
The exact tail asymptotics of the span $\chi_2(\xi)$ when $\xi$ is a
centered stationary Gaussian process with covariance function that
satisfies certain regular conditions is obtained in \cite{PitPris81}.
The same result should be retrieved, using first a time scaling
and then resorting to Remark \ref{remarkmain}. This observation is
confirmed in {Proposition~\ref{PP}} below.

Hereafter, assume that $\{\xi(t),t\in[0,1]\}$ is a centered stationary
Gaussian process with covariance function $r_\xi(s)$ satisfying the
following conditions:
\begin{longlist}[S1$'$:]
\item[S1$'$:] $r_\xi(t)$ attains its minimum on $[0,1]$ at unique
point $t_m\in(0,1)$;

\item[S2$'$:] there exist positive constants $a_1, a_2,\alpha
_1$ and
$\alpha_2\in(0,2)$ such that
\[
r_\xi(t)=r_\xi(t_m)+a_1\llvert
t-t_m \rrvert ^{\alpha
_1}\bigl(1+o(1)\bigr),\qquad t\rightarrow
t_m
\]
and
\[
r_\xi(t)=1-a_2t^{\alpha_2}\bigl(1+o(1)\bigr),\qquad t
\rightarrow0; %
\]
\item[S3$'$:] $r_\xi(t)<1$ for any $t\in(0,1]$.
\end{longlist}

\begin{prop}\label{PP} Let $\{\xi(t),t\in[0,1]\}$ be a centered
stationary Gaussian process with covariance
function $r_\xi(t)$ satisfying \textup{S1$'$--S3$'$}.
If $r_\xi(t)$ is twice continuously differentiable on interval
$[t_m-\mu
,t_m+\mu]$ for some positive small constant~$\mu$, 
then, as $u\rightarrow\infty$,
%
\begin{eqnarray}
\label{eq:PP} \mathbb{P} \bigl( \chi_2(\xi)> u \bigr) &=&2\mathbb{P}
\bigl( \chi_1(\xi)> u \bigr)
\nonumber
\\
&=& 2^{2-{4}/{\alpha_2}+{2}/{\alpha_1}} (1-t_m) \mathcal {H}_{\alpha_2}^{2}
a_2^{{2}/{\alpha_2}} \bigl(1-r_\xi(t_m)
\bigr)^{2-{4}/{\alpha
_2}+
{2}/{\alpha_1} }
\\
&&{}\times u^{{4}/{\alpha_2}-{2}/{\alpha_1}} \Psi \biggl(\frac{u}{\sqrt{2(1-r_\xi(t_m))}} \biggr) \bigl(1+o(1)
\bigr).\nonumber
\end{eqnarray}
\end{prop}

\begin{prop}\label{BB2} If $\{X(t),t\in[0,1]\}$ is the Brownian bridge
given in~\eqref{eq:Z2}, then, as $u\rightarrow\infty$,
%
\begin{equation}
\label{eq:BB2} \mathbb{P} \bigl( \chi_2(X)> u \bigr) = {2}\mathbb{P}
\bigl( \chi_1(X)> u \bigr) =2^{{9}/{2}}\sqrt{\pi }u^3
\Psi (2u) \bigl(1+o(1)\bigr).
\end{equation}
\end{prop}

\begin{remarks}
(a) The claim in \eqref{eq:PP} is consistent with Theorem~2.1 in
\cite{PitPris81}.

(b) Let $B_\alpha$ be a standard fBm and consider its maximum loss
$\chi
_1(B_\alpha)$ and span $\chi_2(B_\alpha)$.
The variance function of the random field $X_1(s,t):=B_\alpha
(t)-B_\alpha(s)$
is given by
\[
\sigma_{X_1}^2(s,t)=\llvert t-s \rrvert ^{\alpha},\qquad
(s,t)\in[0,1]^2
\]
and
attains its maximum only at points $(0,1)$ and $(1,0)$.
Therefore, Theorem~8.2 in \cite{Pit96} yields that, as $u\rightarrow
\infty$,
\begin{longlist}[(iii)]
\item[(i)] if $\alpha\in(0,1)$
\[
\mathbb{P} \bigl( \chi_2(B_\alpha)> u \bigr) =2\mathbb{P}
\bigl( \chi_1(B_\alpha)> u \bigr) = 2^{3-2/\alpha
}\alpha
^{-2}\mathcal{H}_{\alpha}^2 u^{{4}/{\alpha}-4} \Psi(u)
\bigl(1+o(1)\bigr);
\]
\item[(ii)] if $\alpha=1$
\[
\mathbb{P} \bigl( \chi_2(B_\alpha)> u \bigr) =2\mathbb{P}
\bigl( \chi_1(B_\alpha)> u \bigr) = 8 
\Psi(u) \bigl(1+o(1)\bigr);
\]
\item[(iii)] if $\alpha\in(1,2)$
\[
\mathbb{P} \bigl( \chi_2(B_\alpha)> u \bigr) =2\mathbb{P}
\bigl( \chi_1(B_\alpha)> u \bigr) =2 
\Psi(u) \bigl(1+o(1)\bigr).
\]
\end{longlist}
\end{remarks}

\section{Proofs}\label{sec5}
\mbox{}
\begin{pf*}{Proof of Lemma \ref{Lem2}} The claim follows by a direct
application of {Lemma~\ref{Lem1}} given in the \hyperref[app]{Appendix}.
\end{pf*} 

\begin{pf*}{Proof of Theorem \ref{thM}} As it will be seen at the
end of the proof, by
symmetry, cases (vi) and (vii) follow from the claims of cases (iii)
and (v), respectively. Thus, we shall first focus on the proof
of cases (i)--(v). In view of Assumption \ref{assA1} there exist some
$\theta\in(0,1)$ and $\rho_0\ge\rho_1$ ($\rho_1$ is as in \ref{assA3})
such that
\[
\sup_{(s,t)\in[0,S]\times[0,\rho_0]}\sigma(s,t)<\theta.
\]
For $\delta(u)=(\ln u/u)^{2/\beta}, u>0$, we may write
\begin{eqnarray*}
&&\mathbb{P} \Bigl( \sup_{(s,t)\in[0,S]\times[T-\delta
(u),T]}X(s,t)>u \Bigr) \\
&&\qquad\le\mathbb{P}
\Bigl( \sup_{(s,t)\in
[0,S]\times[0,T]}X(s,t)>u \Bigr)
\\
&&\qquad \le\mathbb{P} \Bigl( \sup_{(s,t)\in[0,S]\times[T-\delta
(u),T]}X(s,t)>u \Bigr) +\pi
_1(u)+\pi_2(u),
\end{eqnarray*}
where
\begin{eqnarray*}
\pi_1(u)&:=&\mathbb{P} \Bigl( \sup_{(s,t)\in[0,S]\times[0,\rho
_0]}X(s,t)>u
\Bigr) ,
\\
\pi_2(u)&:=&\mathbb{P} \Bigl( \sup_{(s,t)\in[0,S]\times[\rho
_0,T-\delta(u)]}X(s,t)>u
\Bigr).
\end{eqnarray*}
We shall mainly focus on the analysis of 
\begin{equation}
\label{pi} \pi(u):=\mathbb{P} \Bigl( \sup_{(s,t)\in[0,S]\times[T-\delta
(u),T]}X(s,t)>u
\Bigr) ,\qquad u\rightarrow\infty
\end{equation}
and show that for $i=1,2$
%
\begin{equation}
\label{eq:pi12} \pi_i(u)=o\bigl(\pi(u)\bigr),\qquad u\rightarrow\infty,
\end{equation}
which then implies 
\[
\mathbb{P} \Bigl( \sup_{(s,t)\in[0,S]\times[0,T]}X(s,t)>u \Bigr) =\pi(u)
\bigl(1+o(1)\bigr),\qquad u\rightarrow\infty.
\]

The asymptotics of \eqref{pi} will be investigated for the cases
(i)--(v) separately by using a case-specific approach.

\textit{Case} (i) $\beta>\max(\alpha_1,\alpha_2)$: 
For space saving, we consider only the case that $\alpha_1=\alpha
_2=:\alpha$; the other cases can be shown
with similar arguments.
Following the idea of \cite{Pit2001} choose first a constant $\alpha
_0\in(\alpha,\beta)$, and denote
\[
\triangle_{ij}= {\triangle}_i \times{
\triangle}_j,\qquad \triangle _{ij}^T= {
\triangle}_i \times(T-{\triangle}_j)
\]
with
\[
{\triangle}_i=\bigl[iu^{-{2}/{\alpha_0}}, (i+1)u^{-{2}/{\alpha
_0}}
\bigr],\qquad i=0, 1,\ldots.
\]
Set further
\[
\tilde{N}_1(u)= \bigl\lfloor Su^{{2}/{\alpha_0}} \bigr\rfloor +1,\qquad
\tilde{N} _2(u)= \bigl\lfloor(\ln u)^{{2}/{\beta}}u^{{2}/{\alpha
_0}-
{2}/{\beta}}
\bigr\rfloor+1,
\]
where $\lfloor\cdot\rfloor$ stands for the ceiling function.
By Bonferroni's inequality, we have that
%
\begin{eqnarray}
\label{eq:M1}&& \sum_{i=0}^{\tilde{N}_1(u)}\sum
_{j=0}^{\tilde{N}_2(u)}\mathbb{P} \Bigl( \sup
_{(s,t)\in\triangle_{ij}^T }X(s,t)>u \Bigr)\nonumber\\
&&\qquad \ge \pi(u)
\\
&&\qquad\ge\sum_{i=0}^{\tilde{N}_1(u)-1}\sum
_{j=0}^{\tilde
{N}_2(u)-1}\mathbb{P} \Bigl( \sup
_{(s,t)\in\triangle
_{ij}^T}X(s,t)>u \Bigr) -\Sigma_1(u),\nonumber
\end{eqnarray}
with
\begin{eqnarray*}
&&\Sigma_1(u)= \mathop{\ssum_{0\le i,i'\le\tilde{N}
_1(u)-1, 0\le j, j'\le\tilde{N}_2(u)-1}}_{(i,j)\neq(i',j')}
\mathbb {P} \Bigl( \sup_{(s,t)\in\triangle_{ij}^T}X(s,t)>u, \\
&&\hspace*{180pt}\sup
_{(s,t)\in
\triangle_{i'j'}^T}X(s,t)>u \Bigr).
\end{eqnarray*}

For any $\varepsilon\in(0,1)$ and all $u$ large [set $b_{\pm
\varepsilon}:=b(1 \pm
\varepsilon)$]
\begin{eqnarray*}
\mathbb{P} \Bigl( \sup_{(s,t)\in\triangle_{ij}^T}X(s,t)>u \Bigr)& \le&\mathbb{P}
\biggl( \sup_{(s,t)\in\triangle_{ij} }\frac
{X(s,T-t)}{\sigma(s,T-t)}>u_{j-} \biggr)
,
\\
\mathbb{P} \Bigl( \sup_{(s,t)\in\triangle_{ij}^T}X(s,t)>u \Bigr) &\ge&\mathbb{P}
\biggl( \sup_{(s,t)\in\triangle_{ij} }\frac
{X(s,T-t)}{\sigma(s,T-t)}>u_{j+} \biggr)
,
\end{eqnarray*}
with
\[
u_{j-}= u\bigl(1+b_{- \varepsilon}\bigl(ju^{-{2}/{\alpha_0}}
\bigr)^{\beta}\bigr),\qquad u_{j+}= u\bigl(1+b_{+ \varepsilon}
\bigl((j+1)u^{-{2}/{\alpha_0}}\bigr)^{\beta}\bigr).
\]
Let $\{\eta_{\pm\varepsilon}(s,t), (s,t)\in[0,\infty)^2\}$ with
$\varepsilon$ as
above be
centered stationary Gaussian random fields with covariance functions
\[
r_{\eta_{\pm\varepsilon}}(s,t)=\exp \bigl(-(1\pm\varepsilon )^{\alpha} \bigl(
\llvert a_1s \rrvert ^{\alpha}+\llvert a_2
t+a_3s \rrvert ^{\alpha} \bigr) \bigr),\qquad (s,t)\in [0,
\infty)^2,
\]
respectively. By Slepian's lemma (see, e.g., \cite{Berman92} or \cite
{AZI}) for all $u$ large
\[
\mathbb{P} \biggl( \sup_{(s,t)\in\triangle_{ij}}\frac
{X(s,T-t)}{\sigma(s,T-t)}>u_{j-}
\biggr) \le%
\mathbb{P} \Bigl( \sup
_{(s,t)\in\triangle_{ij}}\eta _{+\varepsilon}(s,T-t)>u_{j-} \Bigr).
\]
In view of Theorem~7.2 in \cite{Pit96}, as $u\rightarrow\infty$,
%
\begin{eqnarray}
\label{eq:1upper} 
\pi(u) &\le&\sum
_{i=0}^{\tilde{N}_1(u)}\sum_{j=0}^{\tilde
{N}_2(u)}
\mathbb{P} \Bigl( \sup_{(s,t)\in\triangle_{ij}}\eta _{+\varepsilon}(s,T-t)>u_{j-}
\Bigr)
\nonumber
\\
&=&(1+\varepsilon)^2a_1a_2
\mathcal{H}_{\alpha}^2u^{-
{4}/{\alpha_0}} \sum
_{i=0}^{\tilde{N}
_1(u)}\sum_{j=0}^{\tilde{N}_2(u)}
u_{j-}^{{4}/{\alpha}}\Psi (u_{j-}) \bigl(1+o(1)\bigr)
\nonumber
\\
&=&(1+\varepsilon)^2a_1a_2
\mathcal{H}_{\alpha}^2Su^{-
{2}/{\alpha_0}+
{4}/{\alpha
}}\Psi(u)\sum
_{j=0}^{\tilde{N}_2(u)} \exp \bigl(- b_{- \varepsilon
}
\bigl(ju^{{2}/{\beta
}-
{2}/{\alpha_0}}\bigr)^{\beta} \bigr) \\
&&{}\times\bigl(1+o(1)\bigr)
\nonumber
\\
&=&(1+\varepsilon)^2 a_1a_2
\mathcal{H}_{\alpha}^2S u^{
{4}/{\alpha}-
{2}/{\beta
}}\Psi(u)\int
_{0}^\infty\exp \bigl(-b_{- \varepsilon}x^{\beta
}
\bigr)\,dx \bigl(1+o(1)\bigr).\nonumber
\end{eqnarray}
Similarly, we obtain
%
\begin{eqnarray}
\label{e:1lower} &&\sum_{i=0}^{\tilde{N}_1(u)-1}\sum
_{j=0}^{\tilde{N}_2(u)-1}\mathbb {P} \Bigl( \sup
_{(s,t)\in\triangle_{ij}^T }X(s,t)>u \Bigr)
\nonumber
\\
&& \qquad\ge \sum_{i=0}^{\tilde{N}_1(u)-1}\sum
_{j=0}^{\tilde{N}_2(u)-1}\mathbb {P} \Bigl( \sup
_{(s,t)\in\triangle_{ij} }\eta_{-\varepsilon
}(s,T-t)>u_{j+} \Bigr)
\\
&&\qquad \ge(1-\varepsilon)^2a_1a_2
\mathcal{H}_{\alpha}^2S u^{{4}/{\alpha
}-{2}/{\beta}}\Psi(u)\int
_{0}^\infty\exp \bigl(-b_{+
\varepsilon}x^{\beta}
\bigr)\,dx \bigl(1+o(1)\bigr).\nonumber
\end{eqnarray}
Next, we deal with the double sum part $\Sigma_1(u)$.
Denote the distance of two nonempty sets $A,B\subset\mathbb{R}^n$ by
\[
\rho(A,B)=\inf_{x\in A,y\in B} \Vert x-y \Vert ,
\]
with $ \Vert \cdot\Vert  $ the
Euclidean distance.
We see from \eqref{eq:r} that there exists a positive constant
$\rho_3$ such that
%
\begin{eqnarray}
\label{eq:rrho_3} && \tfrac{3}{2} \bigl(\bigl\llvert a_1
\bigl(s-s'\bigr) \bigr\rrvert ^{\alpha
}+\bigl\llvert
a_2\bigl(t-t'\bigr)+a_3
\bigl(s-s'\bigr) \bigr\rrvert ^{\alpha} \bigr)\nonumber \\
&&\qquad\ge 1-r
\bigl(s,t,s',t'\bigr)
\\
&& \qquad\ge\tfrac{1}{2} \bigl(\bigl\llvert a_1\bigl(s-s'
\bigr) \bigr\rrvert ^{\alpha}+\bigl\llvert a_2
\bigl(t-t'\bigr)+a_3\bigl(s-s'\bigr) \bigr
\rrvert ^{\alpha} \bigr)\nonumber
\end{eqnarray}
for $\llvert  s-s' \rrvert \le2\rho_3, \llvert  T-t \rrvert \le2\rho_3$ and $\llvert  T-t' \rrvert \le
2\rho_3$.
It follows further from \eqref{eq:r1} that 
there exists some $\theta_0\in(0,1)$ such that
\[
\mathop{\sup_{{0\le
i,i'\le\tilde{N}
_1(u)-1, 0\le j, j'\le\tilde{N}_2(u)-1}}}_{\rho({\triangle}_i,
{\triangle}_{i'})>\rho_3} \mathop{\sup
_{(s,t)\in\triangle_{ij}^T}}_{ (s',t')\in
\triangle_{i'j'}^T
}r\bigl(s,t,s',t'
\bigr)<\theta_0.
\]
Next, we divide the double sum part $\Sigma_1(u)$ as follows: 
\[
\Sigma_1(u)= \Sigma_{1,1}(u)+\Sigma_{1,2}(u)+
\Sigma_{1,3}(u),\qquad u\ge0,
\]
where $\Sigma_{1,1}(u)$ is the sum taken on $\rho({\triangle}_i,
{\triangle}
_{i'})>\rho_3$,
$\Sigma_{1,2}(u)$ is the sum taken on $\rho(\triangle_{ij}^T,
\triangle_{i'j'}^T)=0$ and
$\Sigma_{1,3}(u)$ is the sum taken on $u^{-{2}/{\alpha_0}}\le\rho
(\triangle_{ij}^T,
\triangle_{i'j'}^T)$ and $\rho({\triangle}_i, {\triangle}_{i'})\le
\rho_3$. We first
give the
estimation of $\Sigma_{1,1}(u)$.
For $\xi(s,t,s',t'):=\overline{X}(s,t)+\overline{X}(s',t')$ we have
%
\begin{equation}
\label{eq:xir} \mathbb{E} \bigl(\xi^2\bigl(s,t,s',t'
\bigr) \bigr)=4-2\bigl(1-r\bigl(s,t,s',t'\bigr)\bigr)
\end{equation}
implying
\[
\mathop{\sup_{0\le i,i'\le\tilde{N} _1(u)-1, 0\le j,
j'\le\tilde{N}_2(u)-1}}_{\rho({\triangle}_i, {\triangle
}_{i'})>\rho_3} \mathop{\sup
_{(s,t)\in\triangle_{ij}^T}}_{ (s',t')\in
\triangle_{i'j'}^T
}\mathbb{E} \bigl(\xi ^2
\bigl(s,t,s',t'\bigr) \bigr)\le4-2(1-
\theta_0)<4.
\]
Further, we have
\begin{eqnarray*}
&&\mathbb{P} \Bigl( \sup_{(s,t)\in\triangle_{ij}^T}X(s,t)>u,\sup
_{(s,t)\in\triangle_{i'j'}^T }X(s,t)>u \Bigr)
\\
&&\qquad\le \mathbb{P} \Bigl( \sup_{(s,t)\in\triangle_{ij}^T}\overline {X}(s,t)>u,\sup
_{(s,t)\in\triangle_{i'j'}^T}\overline{X}(s,t)>u \Bigr)
\\
&&\qquad\le \mathbb{P} \Bigl( \mathop{\sup_{(s,t)\in\triangle_{ij}^T
}}_{(s',t')\in\triangle
_{i'j'}^T}
\xi\bigl(s,t,s',t'\bigr)>2u \Bigr).
\end{eqnarray*}
By Borell--TIS inequality (see \cite{AdlerTaylor} or \cite{Pit96}), for
$u$ sufficiently large
\[
\mathbb{P} \Bigl( \sup_{(s,t)\in\triangle_{ij}^T}X(s,t)>u,\sup
_{(s,t)\in\triangle_{i'j'}^T}X(s,t)>u \Bigr) \le\exp \biggl(-\frac{(u-a)^2}{2-(1-\theta_0)}
\biggr),
\]
where $a=\mathbb{E} (\sup_{(s,t),(s't')\in[0,S]\times[0,T]}
\xi(s,t,s',t') )<\infty
$. 
Thus
%
\begin{equation}
\label{eq:1neg11} \limsup_{u\rightarrow\infty}\frac{\Sigma_{1,1}(u)}{ u^{
{4}/{\alpha
}-
{2}/{\beta}}\Psi(u)}=0.
\end{equation}
The summand of $\Sigma_{1,2}(u)$ is equal to
\begin{eqnarray*}
&&\mathbb{P} \Bigl( \sup_{(s,t)\in\triangle_{ij}^T}X(s,t)>u \Bigr) +\mathbb{P}
\Bigl( \sup_{(s,t)\in\triangle_{i'j'}^T }X(s,t)>u \Bigr)
\\
&&\qquad{}- \mathbb{P} \Bigl( \sup_{(s,t)\in\triangle_{ij}^T\cup\triangle
_{i'j'}^T}X(s,t)>u \Bigr).
\end{eqnarray*}
Since $\rho(\triangle_{ij}^T, \triangle_{i'j'}^T)=0$, we have for
$(s,t)\in\triangle_{ij}^T\cup
\triangle_{i'j'}^T
$ and sufficiently large $u$
\begin{eqnarray*}
u\bigl(1+b_{- \varepsilon}\bigl((j-1)_+u^{-{2}/{\alpha_0}}\bigr)^{\beta
}
\bigr)&=:&\tilde{u}_{j-}\le \frac
{u}{\sigma(s,t)}\le\tilde{u}_{j+}\\
&:=&u
\bigl(1+b_{+ \varepsilon
}\bigl((j+2)u^{-
{2}/{\alpha_0}}\bigr)^{\beta}\bigr).
\end{eqnarray*}
Using again Theorem~7.2 in \cite{Pit96} for the last term, we have
\[
\mathbb{P} \Bigl( \sup_{(s,t)\in\triangle_{ij}^T\cup\triangle
_{i'j'}^T}X(s,t)>u \Bigr) \ge2(1-
\varepsilon )^2a_1a_2 \mathcal{H}
_{\alpha}^2u^{-{4}/{\alpha_0}} \tilde{u}_{j+}^{{4}/{\alpha
}}
\Psi (\tilde{u}_{j+}) \bigl(1+o(1)\bigr)
\]
as $u\rightarrow\infty$. Consequently, noting that for any
$\triangle_{ij}^T$ there
are at most 8 sets of the form $\triangle_{i'j'}^T$ in $[0,S]\times
[T-\delta
(u),T]$ adjacent with it, we conclude that
\begin{eqnarray*}
\Sigma_{1,2}(u)&\le& 8 \sum_{i=0}^{\tilde{N}_1(u)}
\sum_{j=0}^{\tilde{N}
_2(u)} \bigl( 2(1+
\varepsilon)^2a_1a_2 \mathcal{H}_{\alpha}^2u^{-{4}/{\alpha
_0}}
\tilde {u}_{j-}^{{4}/{\alpha}}\Psi(\tilde{u}_{j-})
\\
&&\hspace*{54pt}{}-2(1-\varepsilon)^2a_1a_2
\mathcal{H}_{\alpha}^2u^{-
{4}/{\alpha_0}} \tilde
{u}_{j+}^{{4}/{\alpha}}\Psi(\tilde{u}_{j+}) \bigr)
\bigl(1+o(1)\bigr)
\end{eqnarray*}
and thus similar arguments as in \eqref{eq:1upper} yield
%
\begin{equation}
\label{eq:1neg12} \limsup_{\varepsilon\rightarrow0}\limsup_{u\rightarrow\infty
}
\frac{\Sigma_{1,2}(u)}{
u^{
{4}/{\alpha}-{2}/{\beta}}\Psi(u)}=0. 
\end{equation}
Finally, we estimate $\Sigma_{1,3}(u)$. Since $u^{-2/\alpha_0}\le
\rho
(\triangle_{ij}^T, \triangle_{i'j'}^T)$ and $\rho({\triangle}_i,
{\triangle}_{i'})\le\rho_3$, it
follows in view of \eqref{eq:rrho_3} that
\[
\mathop{\inf_{0\le i,i'\le \tilde{N}_1(u)-1, 0\le j,
j'\le\tilde{N}_2(u)-1}}_{\rho({\triangle}_i,
{\triangle}_{i'})\le\rho
_3} \mathop{\inf
_{(s,t)\in
\triangle_{ij}^T,
(s',t')\in
\triangle_{i'j'}^T}}_ {
{u^{-
{2}/{\alpha
_0}}\le\rho(\triangle_{ij}^T, \triangle_{i'j'}^T)}} \bigl(1-r\bigl(s,t,s',t'
\bigr) \bigr)\ge\frac{1}{2}\nu u^{-
{2\alpha}/{\alpha_0}}
\]
for some positive constant $\nu$, and thus
\[
\mathop{\sup_{0\le
i,i'\le
\tilde{N}_1(u)-1, 0\le j, j'\le\tilde{N}_2(u)-1}}_{\rho({\triangle}_i,
{\triangle}_{i'})\le\rho
_3} \mathop{\sup
_{(s,t)\in
\triangle_{ij}^T, (s',t')\in
\triangle_{i'j'}^T
}}_
{u^{-
{2}/{\alpha
_0}}\le\rho(\triangle_{ij}^T, \triangle_{i'j'}^T)}\mathbb{E} \bigl(\xi ^2
\bigl(s,t,s',t'\bigr) \bigr)\le4-\nu
u^{-
{2\alpha}/{\alpha_0}}.
\]
Consequently, using the Piterbarg inequality (cf.  Theorem~8.1 in \cite
{Pit96} or Theorem~8.1 in \cite{Pit2001}) for the summand of $\Sigma_{1,3}(u)$
we obtain
\begin{eqnarray*}
&&\mathbb{P} \Bigl( \sup_{(s,t)\in\triangle_{ij}^T}X(s,t)>u,\sup
_{(s,t)\in\triangle_{i'j'}^T }X(s,t)>u \Bigr)
\\
&&\qquad\le\mathbb{P} \Bigl( \mathop{\sup
_{(s,t)\in\triangle_{ij}^T}}_{(s',t')\in\triangle_{ij}^T}\xi \bigl(s,t,s',t'
\bigr)>2u \Bigr)
\\
&&\qquad= o \biggl(\exp \biggl(-\frac{1}{16}\nu u^{-2({(\alpha_0-\alpha)}/{\alpha_0})} \biggr)
\biggr) u^{{4}/{\alpha}-{2}/{\beta}}\Psi(u),
\end{eqnarray*}
which implies that
%
\begin{eqnarray}
\label{eq:1neg13} &&\limsup_{u\rightarrow\infty}\frac{\Sigma_{1,3}(u)}{ u^{
{4}/{\alpha
}-
{2}/{\beta}}\Psi(u)}
\nonumber
\\
&&\qquad\le\limsup_{u\rightarrow\infty} \mathop{\ssum_{0\le i,i'\le\tilde
{N}_1(u)-1, 0\le j,
j'\le\tilde{N}_2(u)-1}}_{(i,j)\neq(i',j')} o \biggl(\exp \biggl(-
\frac
{1}{16}\nu u^{-2
({(\alpha_0-\alpha)}/{\alpha_0})} \biggr) \biggr) \\
&&\qquad=0.\nonumber 
\end{eqnarray}
Hence, in view of (\ref{eq:M1})--(\ref{e:1lower}), (\ref
{eq:1neg11})--(\ref
{eq:1neg13}) and by letting $\varepsilon\rightarrow0$ we conclude that
\[
\pi(u)=a_1a_2 \mathcal{H}_{\alpha}^2S
u^{{4}/{\alpha}-
{2}/{\beta
}}\Psi (u)\int_{0}^\infty\exp
\bigl(-b x^{\beta} \bigr)\,dx \bigl(1+o(1)\bigr),\qquad  u\rightarrow\infty.
\]

\textit{Case} (ii) $\beta=\alpha_1=\alpha_2$: 
In order to simplify notation, we set $\alpha
:=\alpha_1=\alpha_2$.
Let $S_1, T_1$ be two positive constants and define
\begin{eqnarray*}
\widehat\Delta_i&=&\bigl[iS_1u^{-{2}/{\alpha}},
(i+1)S_1u^{-
{2}/{\alpha}}\bigr], \qquad i=0,\ldots, N_1(u),
\\
\widetilde\Delta_i&=&\bigl[iT_1u^{-{2}/{\alpha}},
(i+1)T_1u^{-{2}/{\alpha}}\bigr], \qquad {i=0},\ldots, N_2(u),
\\
\overline{\triangle}_{ij}&=& \widehat\Delta_i \times
\widetilde \Delta_j, \qquad\overline{\triangle}_{ij}^T=
\widehat\Delta_i \times (T- \widetilde\Delta_j),
\end{eqnarray*}
where
\[
N_1(u)= \biggl\lfloor\frac{S}{S_1}u^{{2}/{\alpha}} \biggr
\rfloor +1,\qquad  N_2(u)= \biggl\lfloor\frac{(\ln u)^{{2}/{\beta}}}{T_1} \biggr\rfloor+1.
\]

Again, Bonferroni's inequality implies
%
\begin{eqnarray}
\label{eq:Bonf2} &&\Sigma_2(u)+\sum_{i=0}^{N_1(u)}
\mathbb{P} \Bigl( \sup_{(s,t)\in
\overline{\triangle}_{i0}^T}X(s,t)>u \Bigr)\nonumber\\
&&\qquad\ge \pi(u)
\\
&&\qquad\ge\sum_{i=0}^{N_1(u)-1}\mathbb{P} \Bigl(
\sup_{(s,t)\in
\overline{\triangle}_{i0}^T }X(s,t)>u \Bigr) -\Sigma_3(u),\nonumber
\end{eqnarray}
where
\begin{eqnarray*}
\Sigma_2(u)&=&\sum_{i=0}^{N_1(u)}
\sum_{j=1}^{N_2(u)}\mathbb{P} \Bigl( \sup
_{(s,t)\in\overline{\triangle}_{ij}^T}X(s,t)>u \Bigr),
\\
\Sigma_3(u)&=& \mathop{\sum\sum}_
{0\le i<i'\le N_1(u)-1}\mathbb{P}
\Bigl( \sup_{(s,t)\in\overline{\triangle}_{i0}^T}X(s,t)>u,\sup_{(s,t)\in\overline{\triangle}_{i'0}^T}X(s,t)>u
\Bigr).
\end{eqnarray*}
Since our approach is of asymptotic nature, for any fixed $0\le i\le
N_1(u)$, the local structures of the variance and correlation of the
Gaussian random field $X$ on $\overline{\triangle}_{i0}^T$ are the
only necessary properties
influencing the asymptotics. Therefore, 
\begin{eqnarray*}
&&\mathbb{P} \Bigl( \sup_{(s,t)\in\overline{\triangle
}_{i0}^T}X(s,t)>u \Bigr) \\
&&\qquad= \mathbb{P}
\biggl( \sup_{(s,t)\in\overline{\triangle}_{i0}}\frac
{\eta(s,t)}{1+bt^\beta}>u \biggr) \bigl(1+o(1)
\bigr)
\end{eqnarray*}
as $u\rightarrow\infty$, where $\{\eta(s,t), (s,t)\in[0,S]\times
[0,T]\}$ is the
same as in {Lemma~\ref{Lem2}}.
Hence, {Lemma~\ref{Lem2}} implies
%
\begin{equation}
\label{eq:2main} \sum_{i=0}^{N_1(u)}\mathbb{P}
\Bigl( \sup_{(s,t)\in\overline
{\triangle}_{i0}^T}X(s,t)>u \Bigr) = \frac{S}{S_1}u^{{2}/{\alpha}}
\mathcal {H}_{Y_1}^{b}[S_1,T_1]
\Psi(u) \bigl(1+o(1)\bigr)
\end{equation}
as $u\to\infty$.
Similarly,
%
\begin{eqnarray}
\label{eq:3main}&& \sum_{i=0}^{N_1(u)-1}\mathbb{P}
\Bigl( \sup_{(s,t)\in\overline
{\triangle}_{i0}^T}X(s,t)>u \Bigr)
\nonumber
\\[-8pt]
\\[-8pt]
\nonumber
&&\qquad = \frac{S}{S_1}u^{{2}/{\alpha}}
\mathcal {H}_{Y_1}^{b}[S_1,T_1]
\Psi(u) \bigl(1+o(1)\bigr)
\end{eqnarray}
as $u\to\infty$.
Note that, for any $c,d\in\mathbb{R}$
\begin{eqnarray*}
\llvert c+d \rrvert ^{p}&\le&\llvert c \rrvert ^p+
\llvert d \rrvert ^p,\qquad \mbox{if }p\in (0,1],
\\
\llvert c+d \rrvert ^{p}&\le&2^{p-1}\bigl(\llvert c
\rrvert ^p+\llvert d \rrvert ^p\bigr)\qquad \mbox {if } p\in(1,
\infty).
\end{eqnarray*}
In view of Slepian's lemma,
\begin{eqnarray*}
&&\mathbb{P} \Bigl( \sup_{(s,t)\in\overline{\triangle
}_{ij}^T}X(s,t)>u \Bigr)
\\
&&\qquad\le\mathbb{P} \Bigl( \sup_{(s,t)\in\overline{\triangle
}_{ij}}\eta(s,t)>u\bigl(1+b
\bigl(jT_1u^{-{2}/{\alpha}}\bigr)^{\beta}\bigr) \Bigr)
\bigl(1+o(1)\bigr)
\\
&&\qquad\le\mathbb{P} \Bigl( \sup_{(s,t)\in\overline{\triangle
}_{ij}}\tilde{\eta}(s,t)>u\bigl(1+b
\bigl(jT_1u^{-{2}/{\alpha}}\bigr)^{\beta
}\bigr) \Bigr)
\bigl(1+o(1)\bigr)
\end{eqnarray*}
as $u\rightarrow\infty$, where $\{\tilde{\eta}(s,t), (s,t)\in
[0,S]\times
[0,T]\}$
is a centered homogeneous Gaussian random field with covariance function
\[
r_{\tilde{\eta}}(s,t)=\exp \bigl(-\llvert \tilde{a}_1 s \rrvert
^{\alpha
}-\llvert \tilde{a}_2 t \rrvert ^{\alpha} \bigr),\qquad
(s,t)\in[0,S]\times[0,T],
\]
with $\tilde{a}_1=(a_1^{\alpha}+2 \llvert  a_3 \rrvert
^{\alpha})^{1/\alpha}$ and
$\tilde{a}_2=2^{1/\alpha}a_2$.
It follows further, using {Lemma~\ref{Lem2}} that
\begin{eqnarray*}
&&\mathbb{P} \Bigl( \sup_{(s,t)\in\overline{\triangle
}_{ij}^T}X(s,t)>u \Bigr)
\\
&&\qquad\le\mathbb{P} \Bigl( \sup_{(s,t)\in\overline{\triangle
}_{ij}}\tilde{\eta}(s,t)>u\bigl(1+b
\bigl(jT_1u^{-{2}/{\alpha}}\bigr)^{\beta
}\bigr) \Bigr)
\bigl(1+o(1)\bigr)
\\
&&\qquad=\mathcal{H}_{\tilde{Y}_2}^0[S_1,T_1]
\frac{1}{\sqrt{2\pi}u} \exp \biggl(-\frac
{u^2(1+2b(jT_1u^{-{2}/{\alpha}})^{\beta})}{2} \biggr) \bigl(1+o(1)\bigr)
\\
&&\qquad=\mathcal{H}_{\tilde{Y}_2}^0[S_1,T_1]
\exp \bigl(-b(jT_1)^{\beta
} \bigr)\Psi (u) \bigl(1+o(1)\bigr)
\end{eqnarray*}
as $u\to\infty$, where $\mathcal{H}_{\tilde{Y}_2}^0[S_1,T_1]$ is
defined in a
similar way as $\mathcal{H}_{Y_2}^0[S_1,T_1]$ with $a_i, i=1,2$
replaced by
$\tilde{a}_i, i=1,2$.
Consequently, as $u\to\infty$,
%
\begin{equation}
\label{eq:2neg2} \Sigma_2(u)\le\sum_{j=1}^\infty
\frac{S}{S_1}u^{{2}/{\alpha
}}\mathcal{H} _{\tilde
{Y}_2}^0[S_1,T_1]
\exp \bigl(-b(jT_1)^{\beta} \bigr)\Psi(u) \bigl(1+o(1)\bigr).
\end{equation}
From \eqref{eq:r1}, there exists some $\theta_1\in(0,1)$ such that
\[
\mathop{\sup_{1\le i<i'\le
N_1(u)}}_{\rho(\widehat\Delta_i, \widehat\Delta_{i'})>\rho
_3} \mathop{\sup
_{s\in\widehat\Delta_i,
s'\in\widehat\Delta
_{i'}}}_{ t,t'\in[0,T]}r\bigl(s,t,s',t'
\bigr)<\theta_1,
\]
where $\rho_3$ is the same as in \eqref{eq:rrho_3}.
Below we shall re-write $\Sigma_3(u)$ as
\[
\Sigma_3(u) = \Sigma_{3,1}(u)+\Sigma_{3,2}(u)+
\Sigma_{3,3}(u), \qquad u\ge0,
\]
where $\Sigma_{3,1}(u)$ is the sum taken on $\rho(\widehat\Delta
_i, \widehat\Delta
_{i'})>\rho_3$,
$\Sigma_{3,2}(u)$ is the sum taken on $i'=i+1$, and $\Sigma_{3,3}(u)$
is the sum taken on $i'>i+1$ and $\rho(\widehat\Delta_i, \widehat
\Delta_{i'})\le\rho_3$.
First, note that the estimation of $\Sigma_{3,1}(u)$ can be derived
similarly to that of $\Sigma_{1,1}(u)$ in case (a), and thus
for $u$ sufficiently large
%
\begin{equation}
\label{eq:2neg31} 
\Sigma_{3,1}(u)\le\frac{S^2}{S_1^2}
u^{{4}/{\alpha}} \exp \biggl(-\frac{(u-a)^2}{2-(1-\theta_1)} \biggr),
\end{equation}
where $a$ is the same as in \eqref{eq:1neg11}.
Next, we consider $\Sigma_{3,3}(u)$. In view of \eqref{eq:rrho_3} and~\eqref{eq:xir}, it follows that for $s\in\widehat\Delta_i, s'\in
\widehat\Delta_{i'},
t,t'\in T-\widetilde\Delta_0$ and $u$
large enough
%
\begin{equation}
\label{eq:xir2} 2\le\mathbb{E} \bigl(\xi^2\bigl(s,t,s',t'
\bigr) \bigr)\le4-\bigl\llvert a_1\bigl(i'-i
\bigr)S_1 \bigr\rrvert ^{\alpha}u^{-2}.
\end{equation}
Further set $\overline{\xi}(s,t,s',t')=\xi(s,t,s',t')/\sqrt
{\operatorname{Var}
(\xi
(s,t,s',t'))}$.
Following similar argument as in the proof of Lemma~6.3 in
\cite
{Pit96}, we obtain that
\begin{eqnarray*}
&&\mathbb{E} \bigl(\overline{\xi}\bigl(s,t,s',t'
\bigr)-\overline{\xi }\bigl(v,w,v',w'\bigr)
\bigr)^2\\
&&\qquad\le 4\bigl(\mathbb{E} \bigl(\overline{X}(s,t)-
\overline{X}(v,w) \bigr)^2
+\mathbb{E} \bigl(\overline{X}\bigl(s',t'\bigr)-
\overline{X}\bigl(v',w'\bigr) \bigr)^2
\bigr).
\end{eqnarray*}
Moreover, from \eqref{eq:rrho_3} we see that, for $u$ sufficiently large
\[
\mathbb{E} \bigl(\overline{X}(s,t)-\overline{X}(v,w) \bigr)^2
\le3 \bigl(\bigl\llvert \tilde{a}_1(s-v) \bigr\rrvert ^{\alpha}+
\bigl\llvert \tilde{a} _2(t-w) \bigr\rrvert ^{\alpha
} \bigr)
\]
implying thus 
\begin{eqnarray}
\label{eq:xirz} &&\mathbb{E} \bigl(\overline{\xi}\bigl(s,t,s',t'
\bigr)-\overline{\xi }\bigl(v,w,v',w'\bigr)
\bigr)^2
\nonumber
\\[-8pt]
\\[-8pt]
\nonumber
&&\qquad\le2\bigl(1- r_{\zeta}
\bigl(s-v,t-w,s'-v',t'-w'
\bigr)\bigr),
\end{eqnarray}
where
\[
r_{\zeta}\bigl(s,t,s',t'\bigr)=\exp \bigl(-7
\bigl( \llvert \tilde{a}_1s \rrvert ^{\alpha
}+\llvert
\tilde{a}_2t \rrvert ^{\alpha}+\bigl\llvert \tilde
{a}_1s' \bigr\rrvert ^{\alpha}+\bigl\llvert
\tilde{a} _2t' \bigr\rrvert ^{\alpha
} \bigr)
\bigr)
\]
is the covariance function of the homogeneous Gaussian random
field $\{\zeta(s,\break t,s', t')$, $(s,t,s',t')\in(0,\infty)^4\}$. Consequently,
\eqref{eq:xir2}, \eqref{eq:xirz} and Slepian's lemma imply
\begin{eqnarray*}
&&\mathbb{P} \Bigl( \sup_{(s,t)\in\overline{\triangle
}_{i0}^T}X(s,t)>u,\sup
_{(s,t)\in\overline{\triangle}_{i'0}^T
}X(s,t)>u \Bigr)
\\
&&\qquad\le\mathbb{P} \biggl( \mathop{\sup_{(s,t)\in\overline{\triangle}_{i0}^T }}_
{(s',t')\in\overline{\triangle
}_{i'0}^T}
\zeta \bigl(s,t,s',t'\bigr)>\frac{2u}{\sqrt{4-\llvert  a_1(i'-i)S_1 \rrvert
^{\alpha}u^{-2}}}
\biggr).
\end{eqnarray*}
We obtain further from a similar lemma as {Lemma~\ref{Lem2}} (cf.
Lemma~6.1
in \cite{Pit96}) that
\begin{eqnarray*}
&&\mathbb{P} \biggl( \mathop{\sup_{(s,t)\in\overline{\triangle
}_{i0}^T}}_{(s',t')\in\overline{\triangle
}_{i'0}^T}\zeta
\bigl(s,t,s',t'\bigr)>\frac{2u}{\sqrt{4-\llvert  a_1(i'-i)S_1 \rrvert
^{\alpha}u^{-2}}} \biggr)
\\
&&\qquad= \bigl(\tilde{\mathcal{H}}_{\tilde{Y}_2}^0[S_1,T_1]
\bigr)^2\frac{1}{\sqrt
{2\pi}u} \exp \biggl(-\frac{4u^2}{2(4-\llvert  a_1(i'-i)S_1 \rrvert ^{\alpha
}u^{-2})} \biggr)
\bigl(1+o(1)\bigr),
\end{eqnarray*}
where $\tilde{\mathcal{H}}_{\tilde{Y}_2}^0[S_1,T_1]$ is defined in a
similar way
as ${\mathcal{H}}_{ Y_2}^0[S_1,T_1]$ with $ a_1, a_2$ replaced by
$7^{1/\alpha
}\tilde{a}_1, 7^{1/\alpha}\tilde{a}_2$, respectively. Consequently,
for all
large $u$,
%
\begin{equation}\quad
\label{eq:2neg33} \Sigma_{3,3}(u)\le \frac{S}{S_1}\sum
_{j\ge1} \bigl(\tilde{\mathcal{H} }_{\tilde
{Y}_2}^0[S_1,T_1]
\bigr)^2\exp \biggl(-\frac{1}{8}\llvert a_1jS_1
\rrvert ^{\alpha
} \biggr)
u^{{2}/{\alpha}}\Psi(u) \bigl(1+o(1)\bigr).
\end{equation}
Next, we consider $\Sigma_{3,2}(u)$. For any $u$ positive,
\begin{eqnarray*}
&&\mathbb{P} \Bigl( \sup_{(s,t)\in\overline{\triangle
}_{i0}^T}X(s,t)>u,\sup
_{(s,t)\in\overline{\Delta
}_{(i+1)0}^T}X(s,t)>u \Bigr)
\\
&&\qquad\le\mathbb{P} \Bigl( \sup_{(s,t)\in\overline{\triangle
}_{i0}^T}X(s,t)>u,\\
&& \hspace*{44pt}\sup
_{(s,t)\in[(i+1)S_1u^{-{2}/{\alpha}},
(i+1)S_1u^{-{2}/{\alpha}}+\sqrt{S_1}u^{-{2}/{\alpha
}}]\times( T-\widetilde\Delta_0)} X(s,t)>u \Bigr)
\\
&&\qquad\quad{}+\mathbb{P} \Bigl( \sup_{(s,t)\in\overline{\triangle
}_{i0}^T}X(s,t)>u, \\
&&\hspace*{59pt}\sup
_{(s,t)\in[(i+1)S_1u^{-{2}/{\alpha
}}+\sqrt{S_1}u^{-{2}/{\alpha}},(i+2)S_1u^{-{2}/{\alpha
}}]\times( T-\widetilde\Delta_0)} X(s,t)>u \Bigr)
\end{eqnarray*}
and further
\begin{eqnarray*}
&&\mathbb{P} \Bigl( \sup_{(s,t)\in\overline{\triangle
}_{i0}^T}X(s,t)>u,\sup
_{(s,t)\in\overline{\Delta
}_{(i+1)0}^T}X(s,t)>u \Bigr)
\\
&&\qquad\le{\mathcal{H}}_{\tilde{Y}_2}^0[\sqrt{S_1},T_1]
\Psi(u) \bigl(1+o(1)\bigr)
\\
&&\qquad\quad{}+\bigl(\tilde{\mathcal{H}}_{\tilde{Y}_2}^0[
\sqrt{S_1},T_1]\bigr)^2 \exp \biggl(-
\frac
{1}{8}\llvert a_1 \sqrt{S_1} \rrvert
^{\alpha} \biggr)\Psi (u) \bigl(1+o(1)\bigr).
\end{eqnarray*}
Therefore, for all large $u$
%
\begin{eqnarray}
\label{eq:2neg32} &&\Sigma_{3,2}(u)\le \frac{S}{S_1} \biggl({
\mathcal{H}}_{\tilde
{Y}_2}^0[\sqrt {S_1},T_1]
\nonumber
\\[-8pt]
\\[-8pt]
\nonumber
&&\hspace*{66pt}{}+\bigl(\tilde{\mathcal{H}}_{\tilde{Y}_2}^0[
\sqrt{S_1},T_1]\bigr)^2 \exp \biggl(-
\frac{1}{8}\llvert a_1 \sqrt{S_1} \rrvert
^{\alpha} \biggr) \biggr)u^{
{2}/{\alpha
}}\Psi(u) \bigl(1+o(1)\bigr).
\end{eqnarray}
Consequently, from (\ref{eq:Bonf2})--(\ref{eq:2neg31})
and 
(\ref{eq:2neg33})--(\ref{eq:2neg32}), we conclude that for any
$S_i,\break T_i, i=1,2$
\begin{eqnarray*}
 &&{S_1}^{-1}\mathcal{H}_{Y_1}^{b }[S_1,T_1]+
\sum_{j=1}^\infty {S_1}^{-1}
\mathcal{H} _{\tilde{Y}_2}^0[S_1,T_1]\exp
\bigl(-b(jT_1)^{\beta} \bigr)
\\
&&\qquad\ge \limsup_{u\rightarrow\infty}\frac{\pi(u)}{Su^{\alpha
/2}\Psi(u)} \ge\liminf
_{u\rightarrow\infty}\frac{\pi(u)}{Su^{\alpha/2}\Psi
(u)}
\\
&&\qquad\ge S_2^{-1}\mathcal{H}_{Y_1}^{b }[S_2,T_2]-
S_2^{-1} \bigl(\tilde {\mathcal{H} }_{\tilde
{Y}_2}^0[S_2,T_2]
\bigr)^2\sum_{j\ge1}\exp \biggl(-
\frac{1}{8}\llvert a_1jS_2 \rrvert ^{\alpha}
\biggr)
\\
&&\qquad\quad{} - S_2^{-1} \biggl({\mathcal{H}}_{\tilde{Y}_2}^0[
\sqrt {S_2},T_2]+\bigl(\tilde{\mathcal{H}
}_{\tilde{Y}_2}^0[\sqrt{S_2},T_2]
\bigr)^2 \exp \biggl(-\frac{1}{8}\llvert a_1
\sqrt{S_2} \rrvert ^{\alpha} \biggr) \biggr).
\end{eqnarray*}
%
Therefore, by similar arguments as in the proof of Theorem D.2 in
\cite{Pit96}, we conclude that
\[
0<\mathcal{M}^{b}_{Y_1,\alpha_1}\le \limsup_{u\rightarrow\infty}
\frac{\pi(u)}{Su^{\alpha/2}\Psi
(u)}\le\liminf_{u\rightarrow
\infty}\frac{\pi(u)}{Su^{\alpha/2}\Psi(u)}\le\mathcal
{M}^{b}_{Y_1,\alpha_1} <\infty
\]
establishing the claim.

\textit{Case} (iii) $\beta=\alpha_2>\alpha_1$: Note that
$\mathcal{M}
_{Y_2,\beta}^{b}$ can be given in terms of Piterbarg and Pickands
constants as
\[
\mathcal{M}_{Y_2,\beta}^{b}=\lim_{T\rightarrow\infty}\lim
_{
S\rightarrow\infty}\frac
{1}{S} \mathcal{H} _{Y_2}^{b}[S,T]=
a_1a_2 \mathcal{P}_{\alpha_2}^{ba_2^{-\alpha
_2}}
\mathcal{H} _{\alpha_1}.
\]
The proof for this case can be established using step-by-step the same
arguments as in case (ii).

\textit{Case} (iv) $\beta<\alpha_2=\alpha_1$: In order to make use
of the notation introduced
in case~(ii) we set $\alpha:=\alpha_1=\alpha_2$.
First, note that $\delta(u)<T_1u^{-2/\alpha}$, which implies
\begin{eqnarray*}
\pi(u)&\le&\mathbb{P} \Bigl( \sup_{(s,t)\in[0,S]\times
(T-\widetilde\Delta_0)}X(s,t)>u \Bigr)
\\
&\le&\sum_{i=0}^{N_1(u)}\mathbb{P} \Bigl(
\sup_{(s,t)\in\overline
{\triangle}_{i0}^T}X(s,t)>u \Bigr)
\\
&\le&\frac{S}{S_1}u^{{2}/{\alpha}}\mathcal {H}_{Y_1}^{0}[S_1,T_1]
\Psi (u) \bigl(1+o(1)\bigr)
\end{eqnarray*}
as $u\rightarrow\infty$. Further, by Assumptions \ref{assA1} and \ref{assA2} we
have that $\mathbb{E} ((X(s,T))^2 )=1, \forall s\in[0,S]$ and
\[
r\bigl(s,T,s',T\bigr) = 1- \bigl(a_1^{\alpha}+
\llvert a_3 \rrvert ^\alpha \bigr)\bigl\llvert
s-s' \bigr\rrvert ^{\alpha} \bigl(1+ o(1)\bigr)
\]
holds uniformly with respect to $s,s'\in[0,S]$, as $\llvert  s-s'
\rrvert \rightarrow0$.
This means that $\{X(s,T), s\in[0,S]\}$ is a locally stationary
Gaussian process. 
Therefore, in view of Theorem~7.1 in \cite{Pit96},
\begin{eqnarray*}
\pi(u)&\ge&\mathbb{P} \Bigl( \sup_{s\in[0,S] }X(s,T)>u \Bigr)
\\
&=&S\bigl(a_1^{\alpha}+\llvert a_3 \rrvert
^{\alpha}\bigr)^{
{1}/{\alpha}} \mathcal{H} _{\alpha}
u^{{2}/{\alpha}}\Psi(u) \bigl(1+o(1)\bigr), \qquad u\rightarrow\infty.
\end{eqnarray*}
%
Letting $T_1\rightarrow0, S_1\rightarrow\infty$, we conclude that
\[
0< \lim_{u\rightarrow\infty}\frac{\pi(u)}{Su^{{2}/{\alpha
}}\Psi
(u)}=\bigl(a_1^{\alpha}+
\llvert a_3 \rrvert ^{\alpha}\bigr)^{
{1}/{\alpha}} \mathcal{H}
_{\alpha} <\infty.
\]

\textit{Case} (v) $\beta<\alpha_2$ \textit{and} $\alpha_1<\alpha_2$:
The claim follows with identical arguments as in the proof of case
(iv).

In order to complete the proof of cases (i)--(v) we only need to
show \eqref{eq:pi12}, for which it is sufficient to give the following
upper bounds for $\pi_1(u)$ and $ \pi_2(u)$.
By Borell--TIS inequality, for $u$ large enough
%
\begin{equation}
\label{eq:Borell} \pi_1(u) \le \exp \biggl(-\frac{ (u-\mathbb{E} (\sup_{(s,t)\in[0,S]\times[0,\rho_0]}X(s,t) ) )^2}{2\theta
^2}
\biggr).
\end{equation}
%
Further, by Assumption \ref{assA3} applying the Piterbarg inequality we
obtain, as $u\to\infty$
%
\begin{eqnarray}
\label{eq:Piter} \pi_2(u) &\le& \mathcal{Q}u^{{4}/{\gamma}-1}\exp
\biggl(-\frac
{u^2}{2\sigma
^2(T-\delta(u))} \biggr)
\nonumber
\\[-8pt]
\\[-8pt]
\nonumber
&=&\mathcal{Q}u^{{4}/{\gamma}-1}\exp \biggl(-\frac
{u^2}{2} \biggr)\exp
\bigl(-b(\ln u)^2 \bigr) \bigl(1+o(1)\bigr),
\end{eqnarray}
%
where $\mathcal{Q}$ is some positive constant not depending on $u$. Therefore,
the proof
of
cases (i)--(v) is complete.

Next, we consider {cases} (vi)--(vii). We introduce a time
scaling of the Gaussian random field $\{X(s,t),(s,t)\in\mathbf
{E}\}$ by
matrix 
$B=\bigl( {{a_3 \atop 0 } \enskip
{a_2 \atop a_2 }}\bigr)$,
that is, let $Z(s,t):=X((s-t)/a_3,t/a_2)$.
By this time scaling, we have
%
\begin{equation}
\label{piu} \mathbb{P} \Bigl( \sup_{(s,t)\in\mathbf{E}}X(s,t)>u \Bigr) =
\mathbb{P} \Bigl( \sup_{(s,t)\in\mathbf{K}}Z(s,t) >u \Bigr) ,
\end{equation}
where $\mathbf{K}$ is a region on $\mathbb{R}^2$ with vertices at
points $(0,0)$,
$(a_2 T,a_2 T)$, $(a_3S,0$) and
$(a_3S+a_2T, a_2T)$. The Gaussian random field $\{Z(s,t),
(s,t)\in
\mathbf{K}\}$ has the following properties: 

(P1) The standard deviation function $\sigma_Z(s,t)$ of $\{
Z(s,t), (s,t)\in\mathbf{K}\}$ satisfies
\[
\sigma_Z(s,t)=1-\frac{b }{a_2^\beta} (a_2T-t
)^\beta \bigl(1+o(1)\bigr),\qquad t\uparrow a_2T.
\]

(P2) The correlation function $r_Z(s,t,s',t')$ of $\{Z(s,t),
(s,t)\in\mathbf{K}\}$ satisfies
\[
r_Z\bigl(s,t,s',t'\bigr)=1- \biggl(
\bigl\llvert s-s' \bigr\rrvert ^{\alpha
_2}+\biggl\llvert
\frac{a_1}{a_3}\bigl(t-t'\bigr)-\frac{a_1}{a_3}
\bigl(s-s'\bigr) \biggr\rrvert ^{\alpha_1} \biggr) \bigl(1+o(1)
\bigr)
\]
for any $(s,t), (s',t')\in\mathbf{K}$ such that $\llvert  s-s'
\rrvert \rightarrow0$ and $t,
t' \uparrow a_2T$, and further there exists some $\delta_0\in(0,T)$
such that
\[
r\bigl(s,t,s',t'\bigr)<1
\]
holds for any $(s,t), (s',t')\in\mathbf{K}_0$ satisfying $s\neq s'$.
Here, $\mathbf{K}_0$ is a region on $\mathbb{R}^2$ with vertices
at points
$(a_2\delta_0,a_2\delta_0)$, $(a_2T,a_2T)$, $(a_3S +a_2\delta
_0,a_2\delta_0)$ and
$(a_3S+a_2T, a_2T)$.

(P3) There exist positive constants $\mathcal{Q}, \gamma, \rho
_1$ and
$\rho
_2$ such that
\[
\mathbb{E} \bigl(\bigl(Z(s,t)-Z\bigl(s',t'\bigr)
\bigr)^2 \bigr)\le\mathcal{Q}\bigl(\bigl\llvert s-s'
\bigr\rrvert ^\gamma+\bigl\llvert t-t' \bigr\rrvert
^\gamma\bigr)
\]
holds for any $(s,t), (s',t')\in\mathbf{K}$ satisfying
$a_2T-t<\rho_1$,
$a_2T-t'<\rho_1$ and $\llvert  s-s' \rrvert <\rho_2$.

Note that in the above proof the most important structural property of
the set $\mathbf{E}$ is that the segment $\mathbf{L}=\{
(s,t)\in\mathbf{E}\dvtx t=T\}$ is
on the boundary of $\mathbf{E}$, which is also the case for $\{Z(s,t),
(s,t)\in\mathbf{K}\}$. Therefore,
in view of the above properties of $\{Z(s,t), (s,t)\in\mathbf{K}\}
$, the
claims of the cases (vi) and (vii) follow by an application of the
claims of cases (iii)
and (v).
The proof is complete.
\end{pf*}

\begin{pf*}{Proof of Proposition \ref{propSGP}} The variance
function of $Z$ is given by
\[
\sigma^2_Z(s,t)= {2\bigl(1-r_X(t)
\bigr)}
\]
{and} attains its maximum on $[0,S]\times\{T\}$. Therefore, it is
sufficient to consider the asymptotics of
\[
\Pi(u):=\mathbb{P} \Bigl( \sup_{(s,t)\in[0,S] \times[0,
T]}Z^*(s,t)> \tilde{u}
\Bigr) ,\qquad  u\rightarrow\infty, %
\]
with
\[
\tilde{u}:=\frac{u}{ \rho_T}\quad \mbox{and}\quad Z^*(s,t):= \frac
{Z(s,t)}{ \rho_T},
\]
where $\rho_T=\sqrt{2 (1- r_X(T))}>0$. The asymptotics of $\Pi(u)$
follows from {Theorem~\ref{thM}} by checking the Assumptions \ref{assA1}--\ref{assA3}.
The standard deviation function of $Z^*$ satisfies
\begin{eqnarray*}
\sigma_{Z^*}(s,t)&=&\frac{\sqrt{2(1-r_X(t))}}{ \rho_T}
\\
&=&1-\frac{ a_1}{{2(1-r_X(T))}}(T-t)^{\alpha_1}\bigl(1+o(1)\bigr),\qquad t\rightarrow T,
\end{eqnarray*}
whereas for its correlation function we have
%
\begin{eqnarray}
\label{eq:RZ} R_{Z^*}\bigl(s,t,s',t'\bigr)
&=&\frac{r_X(\llvert  s+t-s'-t' \rrvert )-r_X( \llvert
s-s'-t' \rrvert )}{2\sqrt
{(1-r_X(t))(1-r_X(t'))}}
\nonumber
\\[-8pt]
\\[-8pt]
\nonumber
&&{} +\frac{ -r_X(\llvert  s+t-s' \rrvert )+r_X(\llvert  s-s'
\rrvert )}{2\sqrt
{(1-r_X(t))(1-r_X(t'))}}.
\end{eqnarray}
Since $r_X(t)$ is twice continuously differentiable in $[\mu,T]$ and
$\llvert  r^{\prime\prime}_X(T) \rrvert \in(0,\infty)$ for some
constant $\mathcal{Q}_1$, we have
\begin{eqnarray*}
&&\bigl\llvert r_X\bigl(t'\bigr)-r_X
\bigl( \bigl\llvert s-s'-t' \bigr\rrvert
\bigr)+r_X(t)-r_X\bigl(\bigl\llvert
s+t-s' \bigr\rrvert \bigr) \bigr\rrvert
\\
&&\qquad\le\mathcal{Q}_1 \bigl(\bigl\llvert t-t'+s-s'
\bigr\rrvert ^{2}+\bigl\llvert s-s' \bigr\rrvert
^{2}\bigr) \bigl(1+o(1)\bigr)
\end{eqnarray*}
as $t,t'\rightarrow T,\llvert  s-s' \rrvert \rightarrow
0$. Consequently, $\alpha_2\in(0,2)$ implies
%
\begin{eqnarray}
\label{eq:RZ2} &&R_{Z^*}\bigl(s,t,s',t'
\bigr)
\nonumber
\\[-8pt]
\\[-8pt]
\nonumber
&&\qquad =1-
\frac{a_2}{\rho_T^2} \bigl( \bigl\llvert t-t'+s-s' \bigr
\rrvert ^{\alpha_2}+\bigl\llvert s-s' \bigr\rrvert
^{\alpha_2} \bigr) \bigl(1+o(1)\bigr)
\end{eqnarray}
as $t,t'\rightarrow T,\llvert  s-s' \rrvert \rightarrow0$.
Next, for any fixed $\varepsilon_0>0$, we have from S3 that
there exists
some $\theta_0$ such that
\[
r_X\bigl(\bigl\llvert s-s' \bigr\rrvert \bigr)\le
\theta_0<1
\]
for any $s,s'\in[0,S]$ satisfying $\llvert  s-s' \rrvert
>\varepsilon_0$. Further, from
S2 we obtain that there exists some positive constant $\delta_0$
such that
\[
2\sqrt{\bigl(1-r_X(t)\bigr) \bigl(1-r_X
\bigl(t'\bigr)\bigr)}\ge\rho_T^2-
\frac{1-\theta_0}{2}> 0
\]
for any $t,t'\in[\delta_0,T]$. Hence,
%
\begin{equation}
\label{eq:RR} R_{Z^*}\bigl(s,t,s',t'\bigr)
\le\frac{1+\theta_0-2r_X(T)}{\rho_T^2-
{(1-\theta_0)}/{2}}<1
\end{equation}
for any $t,t'\in[\delta_0,T]$, $s,s'\in[0,S]$ satisfying $\llvert  s-s' \rrvert >\varepsilon
_0$, and thus both \ref{assA1} and \ref{assA2} are satisfied. It follows that
\begin{eqnarray*}
\mathbb{E} \bigl(Z^*(s,t)-Z^*\bigl(s',t'\bigr)
\bigr)^2&\le& 2 \mathbb {E} \bigl(\overline{Z}(s,t)-\overline{Z}
\bigl(s',t'\bigr) \bigr)^2
\\
&&{}+ \frac{2}{\rho_T^2} \bigl(\sigma_Z(s,t)-\sigma_Z
\bigl(s',t'\bigr)\bigr)^2.
\end{eqnarray*}
Therefore, the differentiability of $r_X(t)$, assumption S2 and
\eqref{eq:RZ2} imply that there exist some positive constants $\rho
_1,\rho_2, \mathcal{Q}_3,\mathcal{Q}_4$
such that
\begin{eqnarray*}
&&\mathbb{E} \bigl(Z^*(s,t)-Z^*\bigl(s',t'\bigr)
\bigr)^2
\\
&&\qquad\le\mathcal{Q}_3 \bigl(\bigl\llvert t-t'+s-s'
\bigr\rrvert ^{\alpha
_2}+\bigl\llvert s-s' \bigr\rrvert
^{\alpha_2}+ \bigl\llvert t-t' \bigr\rrvert ^{2\min(\alpha_1,1)}
\bigr)
\\
&&\qquad\le\mathcal{Q}_4 \bigl(\bigl\llvert t-t' \bigr
\rrvert ^{\min(2\alpha
_1,\alpha_2)}+\bigl\llvert s-s' \bigr\rrvert
^{\min
(2\alpha_1,\alpha_2)}\bigr)
\end{eqnarray*}
for all $s,s' \in[0,S], t,t'\in[\rho_1, T]$ satisfying $\llvert
s-s' \rrvert <\rho
_2$, hence the proof is complete.\end{pf*} 

\begin{pf*}{Proofs of Propositions \ref{propBB} and \ref{propfBm}}
Note first that
the standard deviation of the incremental random field $Z$ of the
Brownian bridge satisfies 
\begin{equation}
\label{eq:BBVar} \sigma_Z(s,t)=\bigl(t(1-t)\bigr)^{{1}/{2}}=
\tfrac{1}{2}- \bigl(t-\tfrac
{1}{2} \bigr)^2\bigl(1+o(1)
\bigr), \qquad t\rightarrow\tfrac{1}{2}.
\end{equation}
Furthermore, for its correlation function we have
%
\begin{equation}
\label{eq:BBr} r_Z\bigl(s,s',t,t'
\bigr)
=1-2\bigl(\bigl
\llvert t-t'+s-s' \bigr\rrvert +\bigl\llvert
s-s' \bigr\rrvert \bigr) \bigl(1+o(1)\bigr)
\end{equation}
as $t,t'\rightarrow1/2,\llvert  s-s' \rrvert \rightarrow0$.

For the fBm incremental random field $Z$, we have for its standard deviation
\[
\sigma_Z(s,t)=t^{{\alpha}/{2}}=1-\frac{\alpha}{2}(1-t)
\bigl(1+o(1)\bigr),\qquad  t\rightarrow1.
\]
As shown in \cite{Pit2001}, the correlation function $r_Z$ of $Z$ 
satisfies
\[
r_Z\bigl(s,s',t,t'\bigr)
= 1-\tfrac{1}{2}\bigl(\bigl\llvert
t-t'+s-s' \bigr\rrvert ^{\alpha}+\bigl\llvert
s-s' \bigr\rrvert ^{\alpha}\bigr) \bigl(1+o(1)
\bigr)%
\]
as $t,t'\rightarrow1,\llvert  s-s' \rrvert \rightarrow0$.
Hence, for both cases {A1--A3} are fulfilled, and thus the claims
follow by a direct application of {Theorem~\ref{thM}}.
\end{pf*} 

\begin{pf*}{Proofs of Propositions \ref{PP} and \ref{BB2}}
By a linear time change using the matrix $A\in\mathbb{R}^{2 \times
2}$ given by
\[
\label{A} A=\left( \matrix{ 1 & 0
\cr
-1 & 1
} \right)
\]
we have for any $u>0$
\[
\mathbb{P} \bigl( \chi_2(\xi)>u \bigr) =\mathbb{P} \Bigl( \sup
_{(s,t)\in A[0,1]^2}\bigl(\xi(t+s)-\xi(s)\bigr)>u \Bigr).
\]
Here, the set $A[0,1]^2=\{(\tilde s,\tilde t)\dvtx (\tilde s,\tilde t)^\top
=A(s,t)^\top, (s,t)\in[0,1]^2\}$ is bounded and convex. 
The variance function of the random field $\{\xi(t+s)-\xi(s),
(s,t)\in A[0,1]^2\}$ is
$2(1-r_\xi(\llvert  t \rrvert ))$ which attains its unique
maximum on the set
$A[0,1]^2$ on two lines $\mathbf{L}_1=\{(s,t)\in A[0,1]^2\dvtx t=t_m\}$
and $\mathbf{L}_2=\{(s,t)\in A[0,1]^2\dvtx t=-t_m\}$. Note that the
differentiability of $r_\xi(t)$ implies $\alpha_1\ge2>\alpha_2$.
Therefore, the claim in \eqref{eq:PP} follows from Remark \ref
{remarkmain}(b); the conditions therein can be established directly as
in the proof of {Proposition~\ref{propSGP}} except \eqref{eq:r2'} for
$i=1, j=2$,
which can also be confirmed by a similar argument as in \eqref{eq:RR}.
Further, since
\[
\mathbb{P} \bigl( \chi_2(X)>u \bigr) =\mathbb{P} \Bigl( \sup
_{(s,t)\in A[0,1]^2}\bigl(X(t+s)-X(s)\bigr)>u \Bigr)
\]
in view of \eqref{eq:BBVar} and \eqref{eq:BBr} we conclude that the
claim in \eqref{eq:BB2} follows immediately from Remark \ref
{remarkmain}(b), and thus the proof is complete.
\end{pf*} 

\begin{appendix}\label{app}
\section*{Appendix}
Let $ \mathbf{ D} $ be a compact set in $\mathbb{R}^2$ such that
$(0,0)\in\mathbf{ D} $, and let $\{\xi_u(s,t), (s,t)\in\mathbf{
D} \}$, $u>0$
be a family of centered Gaussian random fields with a.s. continuous
sample paths. The next lemma is proved based on the classical approach
rooted in the ideas of \cite{PicandsB,PicandsA} (see
also \cite{dekos14}), Lemma~1; in particular, it
implies the claim
of {Lemma~\ref{Lem2}}.

\begin{lem}\label{Lem1}
Let $d(\cdot)$ be a nonnegative continuous function on $[0,\infty)$ and
let $g(u),u>0$ be a positive function satisfying $\lim_{u\rightarrow
\infty}g(u)/u=1$.
Assume that the variance function $\sigma_{\xi_u}^2$
of $\xi_u $ satisfies the following conditions:
\[
\sigma_{\xi_u}(0,0)=1 \qquad\mbox{for all large }
u,\qquad 
\lim_{u\rightarrow\infty}\sup
_{(s,t)\in\mathbf{ D} } \bigl\llvert u^2\bigl(1-
\sigma_{\xi_u}(s,t)\bigr)-d(t) \bigr\rrvert =0,
\]
and there exist some positive constants $G, \nu, u_0$ such that,
for all $u>u_0$
\[
u^2\operatorname{Var}\bigl(\xi_u(s,t)-
\xi_u\bigl(s',t'\bigr)\bigr)\le G \bigl(
\bigl\llvert s-s' \bigr\rrvert ^\nu+\bigl\llvert
t-t' \bigr\rrvert ^\nu\bigr)
\]
holds uniformly with respect to $(s,t), (s',t')\in\mathbf{ D} $. If further
there exists a centered Gaussian random field $\{Y(s,t), (s,t)\in
(0,\infty
)^2\}$ with a.s. continuous sample paths and $Y(0,0)=0$ such that
\[
\lim_{u\rightarrow\infty}u^2\operatorname{Var}\bigl(
\xi_u(s,t)-\xi _u\bigl(s',t'
\bigr)\bigr)=2\operatorname{Var}\bigl(Y(s,t)-Y\bigl(s',t'
\bigr)\bigr)
\]
holds for all $(s,t), (s',t')\in\mathbf{ D} $, then
%
\begin{equation}
\label{eq:lem1} \mathbb{P} \Bigl( \sup_{(s,t)\in\mathbf{ D} }
\xi_u(s,t)>g(u) \Bigr) =\mathcal{H}_{Y}^d[
\mathbf{ D} ]\Psi \bigl(g(u)\bigr) \bigl(1+o(1)\bigr)
\end{equation}
as $u\to\infty$,
where 
\[
\mathcal{H}_Y^{d}[ \mathbf{ D} ]=\mathbb{E} \Bigl(\exp
\Bigl(\sup_{(s,t)\in\mathbf{ D} } \bigl(\sqrt{2}Y(s,t)-\sigma
_Y^2(s,t)-d(t) \bigr) \Bigr) \Bigr).
\]
\end{lem}
\begin{pf} For large $u$, we have
%
\begin{eqnarray}
\label{eq:lamS1} &&\mathbb{P} \Bigl( \sup_{(s,t)\in\mathbf{ D} }
\xi_u(s,t)>g(u) \Bigr)
\nonumber
\\
&&\qquad =\frac{1}{\sqrt{2\pi}g(u)} \exp \biggl(-\frac{(g(u))^2}{2} \biggr)\int
_{-\infty}^\infty e^{w-
{w^2}/{(2(g(u))^2)}}
\\
&& \qquad\quad{}\times\mathbb{P} \biggl( \sup_{(s,t)\in\mathbf{ D} }\xi
_u(s,t)>g(u) | \xi_u (0,0)=g(u)-\frac{w}{g(u)}
\biggr) \,dw.\nonumber
\end{eqnarray}

Let
\[
R_{\xi_u}\bigl(s,t,s',t'\bigr)=\mathbb{E}
\bigl(\xi_u(s,t)\xi_u\bigl(s',t'
\bigr) \bigr),\qquad (s,t),\bigl(s',t'\bigr)\in \mathbf{ D}
\]
be the covariance function of $\xi_u$.
The conditional random field
\[
\biggl\{\xi_u(s,t) | \xi_u(0,0)=g(u)-
\frac{w}{g(u)},  (s,t)\in \mathbf{ D} \biggr\}
\]
has the same distribution as
\[
\biggl\{\xi_u(s,t)-R_{\xi_u}(s,t,0,0)\xi_u(0,0)+R_{\xi
_u}(s,t,0,0)
\biggl(g(u)-\frac{w}{g(u)} \biggr), (s,t)\in\mathbf{ D} \biggr\}.
\]
Thus, the integrand in \eqref{eq:lamS1} can be rewritten as
\begin{eqnarray*}
&&\mathbb{P} \biggl( \sup_{(s,t)\in\mathbf{ D} } \biggl(
\xi_u(s,t) - R_{\xi_u }(s,t,0,0)\xi_u(0,0)\\
&&\hspace*{44pt}{} +
R_{\xi
_u}(s,t,0,0) \biggl(g(u) - \frac{w}{g(u)} \biggr) \biggr) >
g(u) \biggr)
\\
&&\qquad =\mathbb{P} \Bigl( \sup_{(s,t)\in\mathbf{ D} } \bigl(\chi
_u(s,t) - \bigl(g(u)\bigr)^2\bigl(1 -
R_{\xi_u}(s,t,0,0)\bigr)\\
&&\hspace*{110pt}{} + w\bigl(1 - R_{\xi_u
}(s,t,0,0)\bigr)
\bigr) > w \Bigr) ,
\end{eqnarray*}
where
\[
\chi_u(s,t)=g(u) \bigl(\xi_u(s,t)-R_{\xi_u}(s,t,0,0)
\xi_u(0,0)\bigr).
\]
Next, the following convergence
\[
\bigl(g(u)\bigr)^2\bigl(1-R_{\xi_u}(s,t,0,0)\bigr)-w
\bigl(1-R_{\xi_u}(s,t,0,0)\bigr) \to\sigma _Y^2(s,t)+d(t)
\]
holds as $u\rightarrow\infty$, for any $w\in\mathbb{R}$,
uniformly with respect to
$(s,t)\in\mathbf{ D} $. Moreover,
\begin{eqnarray*}
&&\mathbb{E} \bigl( \bigl( \chi_u(s,t)-\chi_u
\bigl(s',t'\bigr) \bigr)^2 \bigr)
\\
&&\qquad =\bigl(g(u)\bigr)^2 \bigl(\mathbb{E} \bigl( \bigl(
\xi_u(s,t)-\xi _u\bigl(s',t'
\bigr) \bigr)^2 \bigr) - \bigl(R_{\xi_u}(s,t,0,0) -
R_{\xi_u}\bigl(s',t',0,0\bigr)
\bigr)^2 \bigr)
\\
&&\qquad \to2\operatorname{Var}\bigl(Y(s,t)-Y\bigl(s',t'
\bigr)\bigr),\qquad u\rightarrow\infty%
\end{eqnarray*}
holds for any $(s,t),(s',t')\in\mathbf{ D} $. Hence, the claim follows
by using the same arguments as in the proof of Lemma~6.1 in \cite
{Pit96} or those in the proof of Lemma~1 in~\cite{dekos14}.
\end{pf}
\end{appendix}

\section*{Acknowledgements}
We are thankful to the Editor, the Associate Editor and the referees
for their comments and suggestions.

%


%




\printaddresses
\end{document}